\documentclass[11pt]{amsart}
\usepackage{amsmath,amsthm, amscd, amssymb, amsfonts,color,mathtools}
\usepackage{hhline}
\usepackage[mathscr]{eucal}
\usepackage[active]{srcltx}
\usepackage{enumitem}
\usepackage{ucs}
\usepackage{graphicx}
\usepackage{booktabs}
\usepackage{multirow}
\usepackage{epigraph}

\usepackage{color}
\usepackage[active]{srcltx}
\usepackage[all]{xy}

\newtheorem{lema}{Lemma}[section]
\newtheorem{example}[lema]{Example}
\newtheorem{teor}[lema]{Theorem}

\newtheorem{prop}[lema]{Proposition}

\theoremstyle{definition}
\newtheorem{defi}[lema]{Definition}

\newtheorem{obs}[lema]{Remark}
\newtheorem{rmks}[lema]{Remarks}

\newcommand{\V}{{\mathcal{V}_{i,j}}}
\newcommand{\Ss}{{\mathcal{S}}}

\newcommand{\com}{\Delta}
\newcommand{\eps}{\varepsilon}
\newcommand{\gr}{\operatorname{gr}}
\newcommand{\cS}{\mathcal{S}}
\newcommand\co{\operatorname{co}}
\newcommand\op{\operatorname{op}}
\newcommand\cop{\operatorname{cop}}
\newcommand\id{\operatorname{id}}

\newcommand\mo{\operatorname{mod}}
\newcommand\End{\operatorname{End}}

\newcommand{\Z}{{\mathbb Z}}
\newcommand{\N}{{\mathbb N}}
\newcommand{\s}{{\mathbb S}}
\newcommand{\bG}{{\mathbb G}}
\newcommand{\bI}{2\Z+1}
\newcommand{\bP}{2\Z}

\newcommand{\M}{{\mathcal M}}

\newcommand{\ydka}{{}^{K}_{K}\mathcal{YD}}
\newcommand{\yda}{{}^{{A}}_{{A}}\mathcal{YD}}
\newcommand{\ydh}{{}^H_H\mathcal{YD}}
\newcommand{\ydhnm}{{}^{H_{n,m}}_{H_{n,m}}\mathcal{YD}}
\newcommand{\ydhs}{{}^{H^{*}}_{H^{*}}\mathcal{YD}}
\newcommand\toba{{\mathfrak B }}
\newcommand\cal{\mathcal }

\newcommand\Ker{\operatorname{Ker}}

\newcommand{\D}{\operatorname{D}}
\newcommand{\HK}{\operatorname{H}}

\def\hit{\mathfrak{R}}
\newcommand{\J}{{\mathcal J}}

\newcommand{\ku}{\Bbbk}
\newcommand{\eva}[2]{\langle {#1,}{\,#2} \rangle}

\def\pf{\begin{proof}}
\def\epf{\end{proof}}
\def\ot{\otimes}

\begin{document}


\title[Finite-dimensional Nichols algebras over dual Radford algebras]
{Finite-dimensional Nichols algebras over dual Radford algebras}

\author[ D. Bagio, G. A. Garc\'ia, J. M. Jury Giraldi, O. M\'arquez]
{D. Bagio, G. A. Garc\'ia, J. M. Jury Giraldi, O. M\'arquez}

\thanks{
This work was partially supported by
ANPCyT-Foncyt, CONICET}
\address{\noindent D. B. and O. M. : Departamento de Matem\'atica, 
	Universidade Federal de Santa Maria,
	Rio Grande do Sul, Brazil}
	\email{bagio@smail.ufsm.br} 
	\email{oscar.f.marquez-sosa@ufsm.br}
	\address{\noindent J. M. J. G. : 
	Instituto de Matem\'atica e Estat\'istica, 
	Universidade Federal do Rio Grande do Sul, 
	Rio Grande do Sul, Brazil}
	\email{joaomjg@gmail.com}
\address{\noindent G. A. G. : Departamento de Matem\'atica, Facultad de Ciencias Exactas,
Universidad Nacional de La Plata. CONICET. Casilla de Correo 172, (1900)
La Plata, Argentina.}
\email{ggarcia@mate.unlp.edu.ar}


\begin{abstract} 
For $n,m\in \N$, let $H_{n,m}$ be the dual of the Radford algebra of dimension $n^{2}m$.
We present new finite-dimensional Nichols algebras arising from the study of simple Yetter-Drinfeld modules over $H_{n,m}$. 
Along the way, we describe the simple objects in $\ydhnm$ and their projective envelopes. Then, 
we determine those
simple modules that give rise to finite-dimensional Nichols algebras for the case $n=2$. There are 18 possible cases. 
We present by generators and relations the corresponding Nichols algebras on five of these
eighteen cases. As an application, we characterize finite-dimensional Nichols algebras over indecomposable modules
for $n=2=m$ and $n=2$, $m=3$, which recovers some results of the second and third author in the former case, and of 
Xiong in the latter.
\end{abstract}

\maketitle

\begin{small}
\begin{epigraph}{\textit{Cualquier destino, por largo y complicado que sea, 
consta en realidad de un solo momento: 
el momento en que el hombre sabe para siempre qui\'en es.}}{Jorge Luis Borges}
\end{epigraph}
\end{small}

\section{Introduction}

Nichols algebras are algebraic objects that encode certain symmetries in classical and quantum (super) geometric objects. 
They were introduced by Nichols \cite{N} and appeared latter as
quantum analogues of enveloping algebras of the nilpotent part of a Borel subalgebra in the work of Lusztig \cite{L} and Rosso \cite{Ro}. 
Nowadays, they play a crucial role in the classification of (pointed) Hopf algebras \cite{AS2}, see \cite{A} for more details. 

Given a braided vector 
space $(V,c)$, one can construct, in different ways, the Nichols algebra $\toba(V)$ associated with it. 
This construction is functorial and it turns out that the Nichols
algebra is determined by $(V,c)$. If the braiding $c$ is rigid, $\toba(V)$ can be realized as a braided Hopf algebra in the category 
of Yetter-Drinfeld modules $\ydh$ over a Hopf algebra $H$, and via the process of \textit{bosonitation} one gets a new Hopf algebra
$\toba(V)\# H$ given by the smash product and the smash coproduct. This fact is used in the \textit{Lifting Method}
introduced by Andruskiewitsch and Schneider \cite{AS} to classify pointed Hopf algebras. 

Fix $\Bbbk$ an algebraically closed field of characteristic zero. 
Let $H$ be a Hopf algebra, $H_{0}$ its coradical (the sum of all its simple subcoalgebras) and denote by $H_{[0]}$ the 
Hopf coradical, \textit{i.e.} the 
subalgebra of $H$ generated by $H_{0}$. Then, one has the following possibilities: either
$H=H_{0}=H_{[0]}$ and $H$ is cosemisimple, or $H=H_{[0]}\neq H_{0}$ and $H$ is generated by its coradical,
or $H\neq H_{[0]} = H_{0}$ and $H_{0}$ is a proper Hopf subalgebra (\textit{i.e.} $H$ has the Chevalley Property), or 
$H\neq H_{[0]} \neq H_{0}$. For the class of Hopf algebras with the Chevalley property, one has the Lifting Method that was applied
with a lot of success to pointed Hopf algebras over abelian groups and some families of non-abelian groups, see for example \cite{A}, \cite{AS2}, \cite{An},
\cite{AGI}, \cite{GaM} and references therein. For the cosemisimple case and
the case where the coradical generates the Hopf algebra there is no general method. This paper concerns with the last case, where the 
coradical generates a proper Hopf subalgebra. This problem is addressed by Andruskiewitsch and Cuadra in \cite{AC}, where they proposed
a \textit{generalized lifting method}. This method coincides with the lifting method in the case that the coradical is indeed a Hopf subalgebra. 

Let $H$ and $A$ be Hopf algebras. 
We say that $H$ is a \textit{Hopf algebra over} $A$ if
$H_{[0]} \simeq A$ as Hopf algebras.
The 
\textbf{standard filtration} $\{H_{ [n]} \}_{ n\geq 0} $ of $H$ is a filtration defined recursively by
$H_{[n]} = \bigwedge^{n+1}H_{[0]}=H_{[0]}\wedge H_{[n-1]}=\{h\in H:\ \com(h)\in H\ot H_{[0]} + H_{[n-1]}\ot H\}$,
for all $n\geq 0$.
If the coradical $H_{0}$ is a Hopf subalgebra, then  
the coradical filtration coincides with the standard one. 
Assume the antipode $\cS$ of $H$ is injective. Then by \cite[Lemma 1.1]{AC}, it holds
that $H_{[0]}$ is a Hopf subalgebra of $H$ and 
$\{H_{ [n]} \}_{ n\geq 0} $ is a Hopf algebra filtration of $H$. In particular, 
the graded algebra associated with the standard filtration
$\gr H = \bigoplus_{ n\geq 0} H_{ [n]} /H_{ [n-1]}$, with $H_{[-1]}=0$,  is a graded Hopf algebra. Moreover,
the homogeneous projection $\pi : \gr H \to H_{ [0]}$ splits the inclusion
of $H_{ [0]}$ in $\gr H$, the \textit{diagram} $R = (\gr H)^{\co\pi}$ is a Hopf algebra in the
category $_{H_{ [0]}}^{H_{ [0]}}\mathcal{YD}$ of Yetter-Drinfeld modules over $H_{[0]}$ and 
$\gr H \simeq R\#H_{[0]}$ as Hopf algebras. It turns out that 
$R$ is graded and connected. As a consequence of this construction one gets the following theorem:

\smallbreak
\noindent
{\bf Theorem.}
\cite[Theorem 1.3]{AC} 
\emph{Any Hopf algebra with injective antipode is a deformation of
the bosonization of a Hopf algebra generated by a cosemisimple coalgebra by
a connected graded Hopf algebra in the category of Yetter-Drinfeld modules
over the latter.\qed}

\smallbreak
The generalized lifting method is then a 
procedure to describe explicitly any Hopf algebra as above.  
If we are interested in finite-dimensional Hopf algebras, 
the main steps reduce to the following: for a Hopf algebra $A$ generated by its coradical, 
\begin{enumerate}
\item[$(a)$] determine all Yetter-Drinfeld modules $V$ in $\yda$ such that the
Nichols algebra $\toba(V)$ is finite-dimensional;
\item[$(b)$] for such $V$, compute all Hopf algebras $L$ such
that $\gr L\simeq \toba(V)\# A$;
\item[$(c)$] prove that any finite-dimensional 
Hopf algebra over $A$
is generated by the first term of the standard filtration. 
\end{enumerate}

A first explicit example is described in \cite{GJG}, where the authors study Hopf algebras over the smallest 
non-pointed non-semisimple Hopf algebra $H_{2,2}$ without the Chevalley property. 
The complete description of Hopf algebras over  $H_{2,2}$ is then achieved by 
putting together the results of \textit{loc.~cit.}, Xiong \cite{X} --where Nichols algebras over semisimple modules are described,
and Andruskiewitsch-Angiono \cite{AA} where generation in degree one is proved. The complete list can be found in \cite[Theorem 6.87]{X2} or
in \cite[Theorem B]{X}, where one has to drop the assumption on the diagram. 
Following a suggestion of Juan Cuadra, we begin the study of Hopf algebras over the family of Hopf algebras $H_{n,m}$ that are dual 
to Radford algebras $R_{n,m}$. Let $n,m\in \N$ and $\omega$ be a primitive $n$-th root of unity.
Radford algebras are pointed Hopf algebras presented by
$$R_{n,m}:=\Bbbk \langle g,x|\ x^n=1-g^n,\ g^{nm}=1,\ gx=\omega xg\rangle,$$
see Subsection \ref{subsec:Radford-alg}. 
These $n^{2}m$-dimensional Hopf algebras were first 
constructed by Radford \cite[Section 2]{Ra} as
examples of non-commutative and non-cocommutative 
Hopf algebras whose Jacobson radical is not a Hopf ideal.
The duals of these algebras $H_{n,m} = R_{n,m}^{*}$ are described in a more general context in \cite{ACE}, see Subsection \ref{subsec:dualRadford}.
It turns out that $H_{n,m}$ is generated as algebra by its coradical. 

On the goal to describe all finite-dimensional Hopf algebras over $H_{n,m}$, we study the category $\ydhnm$ 
of finite-dimensional Yetter-Drinfeld modules over $H_{n,m}$. This is performed through the study of the category of 
finite-dimensional representations over the Drinfeld double $D(H_{n,m}^{\cop})$. We describe the simple objects, their projective
covers and the composition series of the latter. Then we proceed to determine those simple modules such that 
the corresponding Nichols algebra is finite-dimensional. We do this for the case $n=2$ and $m\geq 2$. 
Here we use the general strategy proposed in \cite{AA}, where finite-dimensional Nichols algebras over basic Hopf 
algebras are related to finite-dimensional Nichols algebras of diagonal type. The main idea is to use braided equivalences
of tensor categories to transport the information to one category to the other. In our case, we have the following equivalences 
$$
_{T_{n,m}}^{T_{n,m}}\mathcal{YD}  
\simeq\ _{R_{n,m}}^{R_{n,m}}\mathcal{YD} \simeq\ _{H_{n,m}}^{H_{n,m}}\mathcal{YD} \simeq\
_{D(H^{\cop}_{n,m})}\mathcal{M},
$$
where $T_{n,m}$ is the generalized Taft algebra. The last equivalence is well-known, the second is given by
the general fact that $_{H}^{H}\mathcal{YD} \simeq\ _{H^{*}}^{H^{*}}\mathcal{YD}$ for any finite-dimensional Hopf algebra 
$H$, and the first follows from the fact that 
$R_{n,m}\simeq (T_{n,m})_{\sigma}$, is a 
$2$-cocycle deformation of $T_{n,m}$. We are able to track down all the important information because we know
explicitly the $2$-cocycle involved in the deformation. 
See Section \ref{sec:D-modules} for more details.
To determine those Nichols algebras over simple modules that are finite-dimensional we use the classification of finite-dimensional
Nichols algebras of diagonal type due to Heckenberger \cite{H2}. We get a list of $18$ possible cases. Out of these, we 
present the Nichols algebras over simple modules by generators and relations in the first $5$ cases. As a consequence, we obtain
the results of \cite{GJG} in Theorem \ref{thm:22}, which corresponds to the case $n=2=m$, 
and complete the results of \cite{HX} in Theorem \ref{thm:23}, which corresponds to the case $n=2$, $m=3$. 
Independently of our work, R. Xiong obtained similar results for the case $n=2$ and $m=p$ a prime number in \cite{X2}. 
In particular, he completed his previous result in \cite{HX} and give a complete list of Hopf algebras over $H_{2,3}$ over
indecomposable objects. The authors became aware of Xiong results after finishing this work.

The paper is organized as follows. In Section \ref{sec:prelim} we recall definitions and facts that are used throughout the paper.
Then, in Section \ref{sec:Radford-alg} we introduce Radford algebras $R_{n,m}$, their duals $H_{n,m}$ and 
describe the Drinfeld doubles $D(H_{n,m}^{\cop})$ by generators and relations. In Section \ref{sec:D-modules} we
describe the simple $D(H_{n,m}^{\cop})$-modules in Theorem \ref{theo:RepSimples}, their projective covers 
in Theorems \ref{teor:cover-n} and \ref{teo:cover}, and the composition series of the latter in Proposition \ref{prop-serie}. Finally,
in Section \ref{sec:Nicholsalg} we study the finite-dimensional Nichols algebras over the simple objects in $\ydhnm$ when $n=2$.
In particular, we determine when these are finite-dimensional in Lemma \ref{lem:hecken-pairs}, see Table \ref{tab:Heck}, give the 
presentation of the first five cases of this table in Propositions \ref{prop:Nichols2.1} -- \ref{prop:Nichols7.2}, 
and finish the paper with an application to the cases
$n=2=m$ in Theorem \ref{thm:22} and $n=2$, $m=3$ in Theorem \ref{thm:23}.

\section*{Acknowledgements}
The authors thank
N. Andruskiewitsch and I. Angiono for 
fruitful discussions on preliminary results of the paper during the \textit{Colloquium on Algebras and Representations -- Quantum 17}, 
held in Talca, Chile on March 2017.
We thank as well Mar\'ia Ronco
 for her warm hospitality during the Colloquium. G. A. G. wants to thank also R. Xiong for sharing his PhD Thesis.

\section{Preliminaries}\label{sec:prelim}
In this subsection, we recall some notions that are used along the paper.

\subsection{Conventions}\label{subsec:conv}
Throughout we work over an algebraically closed field $\Bbbk$ of characteristic zero. 
Given a Hopf algebra $H$ over $\Bbbk$, we write $\com$, $\eps$, $\cS$
for the comultiplication, the counit and the
antipode, respectively. Comultiplication and coactions are written using
the Sweedler notation with summation sign suppressed. 

The set $G(H)=\{h\in H\setminus\{0\}:\ \com(h) = h\ot h\}$
denotes the group
of \emph{group-like elements} of a Hopf algebra $H$.
The coradical $H_0$ of $H$ is
the sum
of all simple subcoalgebras; clearly the 
group algebra $\Bbbk G(H)$ is contained in $H_{0}$. The Hopf algebra 
$H$ is called {\it cosemisimple} if $H = H_0$, and it is called {\it  pointed} if $H_{0}=\Bbbk G(H)$.
For $g,h\in G(H)$, the linear space of $(g,h)${\it -primitive elements} 
 is given by the set
$$
\mathcal{P}_{g,h}(H) :=\{x\in H\mid\com(x)=  x\ot g + h\ot x\}.
$$
The elements on such sets are called simply \textit{skew-primitive elements} if no emphasis on the group-likes $g$ and $h$ is needed.
If $g=1=h$, the linear space $\mathcal{P}(H) = \mathcal{P}_{1,1}(H)$ is called
the set of {\it primitive elements}.


Given a $\Bbbk$-vector space $V$, we denote by 
$\Bbbk\{ v_1,\ldots,v_t \}$ the $\Bbbk$-vector subspace of $V$ spanned by $v_1,\ldots,v_t\in V$. For a $\ku$-algebra $A$,
we write $_{A}\mathcal{M}$ for the category of finite-dimensional left $A$-modules.
Let  $q\in \Bbbk$.  For $n \in \N$, we define
the $q$-numbers by
%
$$  \displaylines{
   {(0)}_q  \, := \;  1 \;\; ,  \quad
    {(n)}_q  \, :=  \;  1 + q + \cdots + q^{n-1}  \; =
    \; {\textstyle \sum\limits_{s=0}^{n-1}} \, q^s,  \cr
   {(n)}_q!  \, := \;  {(0)}_q {(1)}_q \cdots {(n)}_q  := \,  {\textstyle \prod\limits_{s=0}^n} {(s)}_q  \;\; ,  \quad
    {\binom{n}{k}}_{\!q}  \, := \;  \frac{\,{(n)}_q!\,}{\;{(k)}_q! {(n-k)}_q! \,}.} $$

The odd (resp. even) integers number are denoted by $\bI$ (resp. $\bP$). We write 
$\bG_n$ for the set of $n$-th roots of unity and 
$\bG'_n$ for the subset of the primitive ones.

\subsection{Yetter-Drinfeld modules} 
Let $H$ be a Hopf algebra. A \textit{left Yetter-Drinfeld module} $M$ over $H$
is a left $H$-module $(M,\cdot)$ and a left $H$-comodule
$(M,\lambda)$ that satisfies that
%
$$ \lambda(h\cdot m) = h_{(1)}m_{(-1)}\cS(h_{(3)})\ot h_{(2)}\cdot m_{(0)},
\qquad \text{ for all }\ m\in M,\, h\in H. $$
We denote by $\ydh$ the corresponding category. It is a braided monoidal category: for any $M, N \in \ydh$,
the braiding $c_{M,N}:M\ot N \to N\ot M$ is given by
$c_{M,N}(m\ot n) = m_{(-1)}\cdot n \ot m_{(0)}$ for all
$m\in M,n\in N$. A Hopf algebra in $\ydh$ is simply called a \textit{braided} Hopf algebra.

Suppose that $H$ is a finite-dimensional. Then
the category $\ydhs$  is braided equivalent to $\ydh$, see \cite[2.2.1]{AG}.
The tensor equivalence is given as follows:
Let $(h_{i})_{i\in I}$ and $(f_{i})_{i\in I}$ be dual basis of $H$ and $H^{*}$.
If $V \in \ydh$, then $V \in \ydhs$ by
\begin{equation}\label{eq:dualyd}
f\cdot v = f(\cS(v_{(-1)})) v_{(0)}\quad\text{ and }\quad
\lambda(v) = \sum_{i} \cS^{-1}(f_{i})\ot h_{i}\cdot v,\quad \text{ for all }v\in V,\ f\in H^{*}.
\end{equation}

\subsection{Braided Hopf algebras and bosonization} Let $H$ be a Hopf algebra and $B$ a braided Hopf algebra in $\ydh$.
The procedure to obtain a
usual Hopf algebra from $B$ and $H$ is called
the Radford biproduct or Majid \emph{bosonization}, and it is usually
denoted by $B \#H$. It is given by $B \# H = B\otimes H$ as vector spaces, and the
multiplication and comultiplication are given by the smash-product
and smash-coproduct,
respectively. That is, for all $b, c \in B$ and $g,h \in H$ one has that 
\begin{align*}
(b \# g)(c \#h) & = b(g_{(1)}\cdot c)\# g_{(2)}h,\\
\com(b \# g) & =b^{(1)} \# (b^{(2)})_{(-1)}g_{(1)} \ot
(b^{(2)})_{(0)}\# g_{(2)},
\end{align*}
where $\com_{B}(b) = b^{(1)}\ot b^{(2)}$ denotes the comultiplication in $B\in \ydh$.
Clearly, the map $\iota: H \to B\#H$ given by $\iota(h) = 1\#h$ for all
$h\in H$ is an injective Hopf algebra map, and the map $\pi: B\#H \to H$
given by $\pi(b\#h) = \eps_{B}(b)h$ for all $b\in B$, $h\in H$
is a surjective Hopf algebra map such that $\pi \circ \iota = \id_{H} $.
Moreover, it holds that $B = (B\#H)^{\co \pi}
=  \{x\in
B\#H\mid (\pi\ot \id )\com(x)=x\ot 1\}$.
Conversely, let $A$ be a Hopf algebra with bijective antipode and
$\pi: A\to H$ a Hopf algebra epimorphism admitting
a Hopf algebra section $\iota: H\to A$.
Then $B=A^{\co\pi}$ is a braided Hopf algebra in $\ydh$ and $A\simeq B\# H$
as Hopf algebras.

\subsection{Cocycle deformations}
A convolution invertible linear map
$\sigma: H\ot H \to \Bbbk $ is called a
\textit{normalized 2-cocycle} on $H$ if for all $a,b,c \in H$ one has that 
$$ \sigma(b_{(1)},c_{(1)})\sigma(a,b_{(2)}c_{(2)}) =
\sigma(a_{(1)},b_{(1)})\sigma(a_{(2)}b_{(2)},c)  $$
and $\sigma (a,1) = \eps(a) = \sigma(1,a)$,
see \cite[Sec. 7.1]{Mo}.
In particular, the inverse of $\sigma $ is given by
$\sigma^{-1}(a,b) = \sigma(\cS(a),b)$.
Using a $2$-cocycle $\sigma$ it is possible to define
a new algebra structure on $H$, which we denote by $ H_{\sigma} $, by deforming the multiplication. 
Moreover, $H_{\sigma}$ is indeed
a Hopf algebra with
$H = H_{\sigma}$ as coalgebras,
deformed multiplication
$$m_{\sigma}(a,b) = a\cdot_{\sigma}b = \sigma(a_{(1)},b_{(1)})a_{(2)}b_{(2)}
\sigma^{-1}(a_{(3)},b_{(3)})\quad\text{ for all }a,b\in H,$$
and antipode
$\cS_{\sigma}$ given by (see \cite{doi} for details)
$$\cS_{\sigma}(a)=\sigma(a_{(1)},\cS(a_{(2)}))\cS(a_{(3)})
\sigma^{-1}(\cS(a_{(4)}),a_{(5)})\quad\text{ for all }a\in H.$$

\begin{obs}
Let $H$ be a Hopf algebra and $\sigma: H\ot H \to \Bbbk$ a normalized $2$-cocycle on $H$. Then 
by \cite[Theorem 3.7]{MO}, there is an equivalence of braided categories 
$\ _{H_{\sigma}}^{H_{\sigma}}\mathcal{YD}  \simeq \ydh$.
By using different techniques, it is possible to prove that (almost) all known examples of finite-dimensional pointed
Hopf algebras are given by $2$-cocycle deformations of bosonizations of Nichols algebras, see \cite{AGI}, \cite{GM}, \cite{Mk},
\cite{GaM}, \cite{GIV}, \cite{FGM} and references therein. 
Thus, in case the $2$-cocycle is given explicitly,
to study the category $\ydh$ for a finite-dimensional pointed
Hopf algebras, it is enough
to study category of Yetter-Drinfeld modules over the bosonization of a Nichols algebra with $H_{0}$. 
The latter category is \textit{a priori} better understood. Indeed, simple and projective objects
can be described using a generalization of weight modules and Verma modules, see for example
\cite{AnRS}, \cite{AB}, \cite{AAMR}, \cite{RS} and \cite{KR}
for the abelian case, and \cite{PV}, \cite{V} for the non-abelian case.
\end{obs}

\subsection{Nichols algebras}\label{subsec:Nichols}
Let $H$ be a Hopf algebra and let $V \in \ydh$. A braided $\N$-graded 
Hopf algebra $R = \bigoplus_{n\geq 0} R(n) \in \ydh$  is called 
a \textit{Nichols algebra} of $V$ if  $\Bbbk\simeq R(0)$, $V\simeq R(1) \in \ydh$,
 $R(1) = \mathcal{P}(R)$ and 
$R$ is generated as an algebra by $R(1)$.
It holds that a Nichols algebra always exists and it is unique up to isomorphism. It 
 is given by the quotient of the tensor algebra $T (V)$ by the largest 
homogeneous two-sided ideal  and coideal $\mathcal{J}$ satisfying that $\mathcal{J}= \bigoplus_{n\geq 2} \mathcal{J}_{n}(V)$ is generated by 
homogeneous elements of degree $\geq 2$.
It is usually denoted by 
$\toba(V) = \bigoplus_{n\geq 0} \toba^{n}(V) $, where $\toba^{0}(V)=\Bbbk$, $\toba^{1}(V)=V$ and  $\toba^{n}(V) = T^{n}(V)/\mathcal{J}_{n}(V)$ for all 
$n\geq 2$. See  \cite{A}, \cite{AS} for more details. 

Given a 
braided vector space $(V,c)$, one may construct the 
Nichols algebra $\toba(V,c)$ in a way similar to the construction above,
by taking a quotient of the tensor algebra $T(V)$ by the homogeneous 
two-sided ideal given by the kernel of an homogeneous symmetrizer.
Let $\mathbb{B}_{n}$ be the braid
group of $n$ letters.
Since $c$ satisfies the braid
equation, it induces a representation $\rho_{n}:\mathbb{B}_{n} \to GL(T^{n}(V))$
for each $n \geq 2$.  Let $M:\s_{n} \to \mathbb{B}_{n}$ be the Matsumoto section
corresponding to the canonical projection 
$\mathbb{B}_{n}\twoheadrightarrow \s_{n}$, and 
consider the morphisms
$$Q_{n} = \sum_{\sigma \in \s_{n}}\rho_{n}(M (\sigma)) \in \End(T^{n}(V)).$$
Then, the Nichols algebra $\toba(V,c)$ 
is the quotient of $T (V)$ by the two-sided ideal
$\mathcal{J}=\bigoplus_{n\geq 2 }\mathcal{J}_{n}(V)$, with $\mathcal{J}_{n}(V)=\Ker Q_{n}$ for all $n\geq 2$.
One usually simply writes 
 $\toba(V)=\toba(V,c)$ if the braiding $c$ is clear from the context, as in the case of $V\in \ydh$.

There is a useful criterion with skew-derivations to find relations in $\toba(V)$ with $V\in \ydh$:
for $f\in V^{*}$ define $\partial_{f} \in \End T(V)$ by 
$$ \partial_{f}(1)=0,\qquad \partial_{f}(v) = f(v)\ \text{ for all }\ v\in V\quad\text{ and }\quad 
 \partial_{f}(xy)= x\partial_{f}(y) + \sum_{j} \partial_{f_{j}}(x)y_{j}, $$
 where  $c^{-1}(y\ot f) = \sum_{j} f_{j}\ot y_{j}$ and $y \in T(V)$.
 Let $x\in T^{n}(V)$ with $n\geq 2$. If $\partial_{f}(x)=0$ for all $f\in V^{*}$, then $x\in \mathcal{J}_{n}(V)$. 

A \textit{pre-Nichols algebra} $\widetilde{\toba}(V,c)$ of a braided vector space $(V, c)$, as defined by Masuoka, 
is any graded braided Hopf algebra intermediate between $T(V )$ and $\toba(V,c)$, that is any
braided Hopf algebra of the form $T(V )/I$ where $I \subseteq \mathcal{J}(V)$ is a homogeneous Hopf ideal.
In particular,  $\widetilde{\toba}(V,c)$ is a graded braided Hopf
algebra $\widetilde{\toba} (V,c)= \bigoplus_{n\geq 0}\widetilde{\toba}^{n}(V,c)$
such that $\widetilde{\toba}^{0}(V,c) = \Bbbk$, $\widetilde{\toba}^{1}(V,c) = V$ and $\widetilde{\toba}(V,c)$ is generated as an
algebra by $V$.


\section{Radford algebras, duals and Drinfeld doubles}\label{sec:Radford-alg}
In this section we recall the definition of Radford algebras, introduce their duals, certain Drinfeld doubles, 
and prove some formulas that are needed in the sequel. 

\subsection{Radford algebras}\label{subsec:Radford-alg}
Let $n,m\in \N$, and $\omega\in \Bbbk$ be a primitive $n$-th root of unity. 
Write $C_{nm}$ for the cyclic group of order $nm$ with generator $g$. 

\begin{defi}\cite[Section 2.1]{Ra}\label{def:Radfordalg}
The Radford algebra $R_{n,m}$ is the associative algebra generated by the elements $g,x$ subject to 
the following relations 
\begin{equation}\label{rel:radford-algebra}
x^n=1-g^n,\qquad g^{nm}=1,\qquad gx=\omega xg.
 \end{equation}
\end{defi}
It is a Hopf algebra with its coalgebra structure determined by $g$ being group-like and $x$ a $(1,g)$-primitive element, 
\textit{i.~e.}, $\Delta(x)=x\otimes 1+g\otimes x$. 
Consequently, $\varepsilon(g)=1$, $\varepsilon(x)=0$, 
$\cS(g)=g^{-1}$ and $\cS(x)=-g^{-1}x$. 
A linear basis is given by the set
$\{x^jg^i\,:\,0\leq j< n, \,\,0\leq i<nm \}$; in particular, $\dim R_{n,m}= n^{2}m$.

The coradical is  $(R_{n,m})_{0} = \Bbbk C_{nm} $, thus $R_{n,m}$ is pointed for all 
$m,n\in \N$. The graded
Hopf algebra associated with the coradical filtration satisfies that $\gr R_{n,m} \simeq \toba(V) \# \Bbbk C_{nm}$,
where $\toba(V)=\Bbbk [x]/(x^{n})$ is the truncated polynonial algebra, and
$\#$ denotes the bosonization of $\toba(V)$ and $\Bbbk C_{nm}$. Here 
$V = \Bbbk \{x\}$ is a braided vector space called \textit{quantum line}. 
Precisely, $V\in \ _{C_{nm}}^{C_{nm}}\mathcal{YD}$ via $g\cdot x=\omega x$ and $\delta(x)=g\otimes x$.
The Hopf algebra $\toba(V) \# \Bbbk C_{nm}$ is better known as the \textit{generalized Taft algebra}; as an algebra,
it can be presented by 
$$
\toba(V) \# \Bbbk C_{nm} = T_{m,n} =\Bbbk\langle g,x:\  x^{n}=0,\ g^{nm}=1,\ gx=\omega xg\rangle.
$$
If $m=1$, then $T_{1,n} = T_{n}$ is the Taft (Hopf) algebra of dimension $n^{2}$.
Following this terminology, Radford algebras
are \textit{liftings of quantum lines}: they are filtered algebras whose associated graded algebra 
is isomorphic to the bosonization $\toba(V) \# \Bbbk C_{nm}$. In particular, 
they are $2$-cocycle deformations of generalized Taft algebras,  see \cite{GM}, \cite{Mk}, \cite{AGI}. 
Thanks to \cite[Section 4.4]{GM}, it is possible to describe explicitly the $2$-cocycle that 
give the deformation.
By \cite[Proposition 4.2]{GM} we have that 
$R_{n,m} = (T_{n,m})_{\sigma}$, where the $2$-cocycle $\sigma$ is  
the exponentiation 
$\sigma = e^{\eta}$ of a Hochschild $2$-cocycle $\eta$  on $T_{n,m}$ given by
\begin{equation}\label{eq:formula-cocycle}
 \eta(x^{a}g^{r}, x^{b}g^{s}) = \begin{cases}
                                 \begin{array}{cc}
                                  \omega^{br} & \text{if }a+b=n, \\
                                  0 & \text{otherwise.}
                                 \end{array} 
                                \end{cases}
\end{equation}
Note that since $\eta*\eta=0$, one has that  $\sigma= e^\eta=\eps\ot\eps+\eta$ 
and $\sigma^{-1}= e^{-\eta}=\eps\ot\eps-\eta$.

\subsection{Duals of Radford algebras} \label{subsec:dualRadford}
Duals of liftings of quantum lines were explicitly calculated in \cite[Section 3.2]{ACE}. 
We describe below the duals of the Radford algebras $R_{n,m}$, which constitute a particular case.

Fix a primitive $nm$-th root of unity $\xi\in \ku^{\times}$ such that $\xi^{m}=\omega$, and consider the elements $U, X$ and $A$ 
of $R_{n,m}^{\ast}$ defined by
the formulas 
\begin{align} \label{pairing-ACE}
&\eva{U}{x^jg^t}=\omega^t\delta_{j,0},& &\eva{X}{x^jg^t}=\delta_{j,1},& &\eva{A}{x^jg^t}=\xi^{t}\delta_{j,0}.&
\end{align}

\begin{prop}\cite[Proposition 3.3]{ACE}
The Hopf algebra $H_{n,m}= R_{n,m}^{*}$ is generated by the elements $U$, $X$ and $A$ 
satisfying the following relations 
\begin{align}
\label{rel:radforddual-algebra1}  &U^n=1,&  &X^n=0,& &A^m=U,\\
\label{rel:radforddual-algebra2}  &UX=\omega XU,&   &UA=AU, & &AX=\xi XA.
\end{align}
Its coalgebra structure is determined by  $U\in G(H_{n,m})$, $X\in \mathcal{P}_{1,U}(H_{n,m})$ and
\[ \com(A)=A\otimes A+\sum\limits_{k=1}^{n-1}\gamma_{n,k}X^{n-k}U^kA\otimes X^kA,\,\,\, 
\text {where } \gamma_{n,k}=\frac{1-\xi^n}{(k)!_{\omega}(n-k)!_{\omega}}.\]
The counit and the antipode of $A$ are defined respectively by $\varepsilon(A)=1$ and $\cS(A)=A^{m-1}U^{n-1}$. 
Note that, as $A\cS(A)=1=\cS(A)A$, one has that $\cS(A)=A^{-1}$. \qed
\end{prop}

\begin{obs}\label{rmk:Hopf-pairing}\rm{
Using the structure given by the proposition above, a straightforward calculation shows that there
exists a Hopf pairing $\eva{-}{-}: H_{n,m}\times R_{n,m} \to \Bbbk $, 
induced by the evaluation, see \cite[Proof of Proposition 3.3]{ACE}. Explicitly, it is given by}
 $$
 \eva{X^{a}A^{b}}{x^c g^d}= \delta_{a,c}\xi^{bd}(a)_{\omega}! \qquad\text{ for all }0\leq a, c<n\text{ and }0\leq b, d <nm.
 $$
 \end{obs}

\begin{obs}\label{rmk:coradical-dual-radford}\rm{
By \cite[Proposition 2.3]{Ra}, $H_{n,m}$ is isomorphic as coalgebra to the direct sum 
$T_{n} \oplus D_{1}\oplus \cdots \oplus D_{m-1}$, where
$D_{i}$  is isomorphic to the simple comatrix coalgebra  $M^{*}(n,\Bbbk)$ of dimension $n^{2}$, for all $1\leq i \leq m-1$. In particular,
$(H_{n,m})_{0} = \Bbbk C_{n} \oplus D_{1}\oplus \cdots \oplus D_{m-1}$ and $\dim (H_{n,m})_{0} = n^{2}(m-1) + n$. As 
$\dim (H_{n,m})_{0}$ does not divide $\dim H_{n,m}$, $(H_{n,m})_{0}$ is not a Hopf subalgebra and consequently $H_{n,m}$ does not have 
the Chevalley Property. Moreover, the subalgebra of $H_{n,m}$ generated
by $(H_{n,m})_{0}$ coincides with $H_{n,m}$. 
For a more general family sharing the same properties, see \cite[Section 3]{ACE}.} 
\end{obs}

\subsection{The Drinfeld double of $H_{n,m}^{\cop}$} Before we introduce our object of study, we first
recall the definition of the Drinfeld double of a finite-dimensional Hopf algebra. 

Let $K$ be a finite-dimensional Hopf algebra.  
There exist a left action of $K$ on its dual $K^{\ast}$, and a right action of $K^{\ast}$ on $K$ which are given respectively by: 
\begin{align*}
& a\rightharpoonup f= \langle f_{(3)}\Ss^{-1}(f_{(1)}), a\rangle f_{(2)},& & a \leftharpoonup f= \langle f, \Ss^{-1}(a_{(3)})a_{(1)}\rangle a_{(2)},&
\end{align*} 
for all $a\in K$ and $f\in K^{\ast}$. 

\begin{defi}\cite[Definition 10.3.5]{Mo}.
The Drinfeld double of $K$ is the Hopf algebra 
$D(K)=(K^{\ast})^{\cop}\bowtie K$, whose underlying $\ku$-vector space is $K^{\ast}\otimes K$, and
the algebra and coalgebra structures are given respectively by: 
\begin{align*} 
&\,\,\,(f\bowtie a)(\tilde{f}\bowtie b) = f(a_{(1)} \rightharpoonup \tilde{f}_{(2)}) 
\bowtie (a_{(2)}\leftharpoonup \tilde{f}_{(1)})b,& &\quad 1 = \varepsilon\bowtie 1,& \\[0.1cm] 
&\,\,\,\Delta(f \bowtie a)=f_{(2)}\bowtie a_{(1)}\otimes f_{(1)}\bowtie a_{(2)},& &\quad\varepsilon(f \bowtie a)=f(1)\varepsilon(a),& 
\end{align*}
for all $f,\tilde{f}\in K^{*}$ and $a,b\in K$.  
\end{defi}

The map from $K$ (resp. $(K^{\ast})^{\cop}$)  to $D(K)$  given by $a\mapsto 1\bowtie a$ (resp. $f\mapsto f \bowtie 1$) 
is a Hopf algebra monomorphism. Thus, $K$ and $(K^{\ast})^{\cop}$ are Hopf subalgebras of $D(K)$. 

The next well-known result establishes the connection between the category $\ydka$ and  the category 
$_{D(K^{\cop})}\cal{M} $ of left $D(K^{\cop})$-modules.

\begin{prop}\cite[Prop. 10.6.16]{Mo}\label{prop:equiv-D-mod-YD}
	There exists an equivalence $ F: \ _{D(K^{\cop})}\cal{M} \to \ydka$ of braided tensor categories. For $V \in\!\!\! \ _{D(K^{\cop})}\cal{M}$, the image 
	$F(V) \in \ydka$ equals $V$ as $K$-module, and the $K$-comodule structure $\lambda:V \to K\ot V $, $v\mapsto v_{(-1)}\otimes v_{(0)}$, is uniquely 
	determined by the equality
	$f\cdot v = f(v_{(-1)})v_{(0)} $ for all $f\in K^{*}$ and $v\in V$.\qed
\end{prop}

\bigbreak
From now on, we set $D:=D(H_{n,m}^{\cop})= (R_{n,m})^{\op\,\cop}\otimes H_{n,m}^{\cop}$. 
In particular, for all $a\in R_{n,m}$ and $f \in H_{n,m}^{\cop}$ we have that 
\begin{align*}
f\rightharpoonup a=\eva{f}{S(a_{(1)})a_{(3)}}a_{(2)}, \qquad f\leftharpoonup a=\eva{S(f_{(1)})f_{(3)}}{a}f_{(2)},\\[2.5pt]
(1\bowtie f)(a\bowtie 1) = (f_{(2)}\rightharpoonup a_{(2)})\bowtie (f_{(1)} \leftharpoonup a_{(1)}) .
\end{align*}

Under the considerations above, 
the next technical lemma is straightforward.
\begin{lema} \label{lem:ident}
	The following identities hold in $D$:
	\begin{align*}
	&U\rightharpoonup g=g,& &X\rightharpoonup g=0,& & A\rightharpoonup g=g,&\\
	&U\leftharpoonup g=U,& &X\leftharpoonup g=\omega^{-1}X,& & A\leftharpoonup g=A,&\\
	& U\rightharpoonup x=\omega^{-1}x,&  & X\rightharpoonup x=\omega^{-1}(g-1),& &A\rightharpoonup x=\xi^{-1}x,& \\
	&U\leftharpoonup x=0,& &X\leftharpoonup x=\omega^{-1}(U-1),& & X^{n-1}UA\leftharpoonup g=\omega X^{n-1}UA,&
	\end{align*}
$$	A\leftharpoonup x=\gamma_{n,1}\xi^{-1}X^{n-1}(U-1)A.$$	
\qed
\end{lema}

The following proposition yields the presentation of $D$ by generators and relations.

\begin{prop}\label{prop:drinfeld_double} The Hopf algebra $D$ is generated by the elements $g,\,x,\,U,\,X$ and $A$ 
satisfying the relations \eqref{rel:radford-algebra} with the opposite product, \eqref{rel:radforddual-algebra1}, 
\eqref{rel:radforddual-algebra2} and the following ones:
	\begin{align}
	\label{eqdouble1}& Ag=gA, \qquad gX= \omega Xg,&\\
	\label{eqdouble2}&xX=\omega Xx+(1-Ug),&  \\
	\label{eqdouble3}&xA=\xi Ax+\gamma_{n,1}X^{n-1} (1-\omega^{n-1}Ug)A.&
	\end{align}
	In particular, $Ug=gU$ and $xU=\omega Ux$.
\end{prop}

\pf 
By the definition of the product in $D$ and Lemma  \ref{lem:ident} we have that
\begin{align*}
 Ag &= (1 \bowtie A)(g\bowtie 1) = (A_{(2)}\rightharpoonup g)\bowtie (A_{(1)}\leftharpoonup g)\\
 & = (A\rightharpoonup g)\bowtie (A\leftharpoonup g)+\sum\limits_{k=1}^{n-1}\gamma_{n,k} (X^kA \rightharpoonup g)\bowtie (X^{n-k} U^kA\leftharpoonup g) \\
 & = g\bowtie A = gA.
\end{align*}
In a similar way, we have
\begin{align*}
Xg &= (1 \bowtie X)(g\bowtie 1) = (X_{(2)}\rightharpoonup g)\bowtie (X_{(1)}\leftharpoonup g)\\
&=(X\rightharpoonup g)\bowtie (U\leftharpoonup g)+(1\rightharpoonup g)\bowtie (X\leftharpoonup g)  = g\bowtie \omega^{-1}X = \omega^{-1}gX. 
 \end{align*}
 Since $U=A^m$, from \eqref{eqdouble1} it follows that $Ug=gU$.  
Analogously, \eqref{rel:radforddual-algebra1} and \eqref{eqdouble3} imply that $xU=\omega Ux$. 
For \eqref{eqdouble2} we have
\begin{align*}
Xx &= (1 \bowtie X)(x\bowtie 1)= (X_{(2)}\rightharpoonup x_{(2)})\bowtie (X_{(1)}\leftharpoonup x_{(1)})\\
&=(X\rightharpoonup x)\bowtie (U\leftharpoonup g)+(1\rightharpoonup x)\bowtie (X\leftharpoonup g)  \\
&\quad+(X\rightharpoonup 1)\bowtie (U\leftharpoonup x)+(1\rightharpoonup 1)\bowtie (X\leftharpoonup x) \\
&=\omega^{-1}(g-1)\bowtie U+\omega^{-1}x\bowtie X+\omega^{-1}1\bowtie (U-1)\\
&=\omega^{-1}(xX+gU-1) = \omega^{-1}(xX+Ug-1).
\end{align*}
The proof of the remaning relations are left to the reader.
 \epf

We end this section with the following equalities that will be used later.

\begin{lema}\label{ref:lema29} For any $k\geq 1$, set $\Gamma(k)= (1- w^k Ug)\in D$. The following equalities hold
\begin{align}
\label{id-gammax} & X^rx^s\Gamma(k)=\Gamma(k-2r+2s)X^rx^s,\,\quad \text{for all }\, r,s\geq 0, &\\
\label{id-xX} &xX^k=\omega^kX^kx+(k)_\omega X^{k-1}\Gamma(k-1), & \\
\label{id-xX1} & x^kX=\omega^kXx^k+(k)_\omega\Gamma(k-1)x^{k-1},&\\
\label{id-xkXkA} & (X^kx^k)A=A(X^kx^k).&
\end{align}
\end{lema}

\pf The identity \eqref{id-gammax} is immediate. We prove \eqref{id-xX} by induction on $k$; the proof of \eqref{id-xX1} 
is similar.
For $k=1$, \eqref{id-xX} is equal to \eqref{eqdouble2}. Let $k\geq 1$ and assume that the formula holds for $k-1$. 
Using that $X^\ell\Gamma(j)=\Gamma(j-2\ell)X^\ell$, for all $j,\ell\geq 1$, we have
\begin{align*}
xX^k&=(\omega^{k-1}X^{k-1}x+(k-1)_{\omega}X^{k-2}\Gamma(k-2))X\\
    &=\omega^{k-1}X^{k-1}(\omega Xx+\Gamma(0)) + (k-1)_{\omega}X^{k-2}X\Gamma(k)\\
    &=\omega^{k}X^{k}x+X^{k-1}(\omega^{k-1}\Gamma(0)+(k-1)_{\omega}\Gamma(k))\\
    &=\omega^kX^kx+(k)_\omega X^{k-1}\Gamma(k-1).
\end{align*}
For \eqref{id-xkXkA} we also proceed by induction on $k$. From \eqref{eqdouble3}, 
$X(xA)=X(\xi Ax)=AXx$  which implies the identity for $k=1$. 
Assume that the result is true for $k>1$. Then 
\begin{align*}
X^{k+1}x^{k+1}A &= X^k(Xx^{k})xA \\
&\overset{\eqref{id-xX1}}{=} X^k\left(\omega^{-k}x^kX-\omega^{-k}(k)_\omega \Gamma(k-1)x^{k-1}\right)xA\\
&\overset{\eqref{id-gammax}}{=}\left( \omega^{-k} X^kx^kXx-\omega^{-k}(k)_\omega  \Gamma(k-1-2k)X^kx^{k}\right) A\\
&\overset{(\ast)}{=}A\left( \omega^{-k}  X^kx^kXx-\omega^{-k}(k)_\omega  \Gamma(k-1-2k)X^kx^{k}\right) \\
&\overset{\eqref{id-gammax}}{=}A\left( \omega^{-k}  X^kx^kXx-\omega^{-k}(k)_\omega X^k \Gamma(k-1)x^{k}\right) \\
&=AX^{k+1}x^{k+1}.
\end{align*}
Here, equality $(\ast)$ follows by induction and the case $k=1$.
\epf


\section{The category $_{D}\mathcal{M}$}\label{sec:D-modules}
In this section we classify the simple objects in $_{D}\mathcal{M}$ and describe their projective covers. 
By \S \ref{subsec:Radford-alg}, we know that $R_{n,m}\simeq (T_{n,m})_{\sigma}$, where $\sigma=e^{\eta}$ is a normalized
$2$-cocycle and $\eta$ is a Hochschild $2$-cocycle as defined in \eqref{eq:formula-cocycle}. Then, 
there is an equivalence of braided categories 
$\ _{T_{n,m}}^{T_{n,m}}\mathcal{YD}  \simeq\  _{R_{n,m}}^{R_{n,m}}\mathcal{YD}$,
given by the identity on the
underlying vector spaces and coactions, but changing the action by the following formula:
$$h \rightharpoonup_{\sigma} v = \sigma (h_{(1)} ,v_{(-1)} )(h_{(2)} 
\rightharpoonup v_{(0)} )_{(0)} 
\sigma^{-1} ((h_{(2)} \rightharpoonup v_{(0)} )_{(-1)} ,h_{(3)} ),$$
for all $h \in R_{n,m},\ v \in V \in\ _{R_{n,m}}^{R_{n,m}}\mathcal{YD}$.
Taking into account the equivalence $\ydh \simeq \ _{H^{*}}^{H^{*}}\mathcal{YD}$ given in \eqref{eq:dualyd} 
and Proposition \ref{prop:equiv-D-mod-YD}, one has
the equivalence of braided tensor categories 
\begin{equation}\label{eq:braided-equiv}
_{D(T^{\cop}_{n,m})}\mathcal{M}\simeq\ _{T_{n,m}}^{T_{n,m}}\mathcal{YD}  
\simeq\ _{R_{n,m}}^{R_{n,m}}\mathcal{YD} \simeq\ _{H_{n,m}}^{H_{n,m}}\mathcal{YD} \simeq\
_{D}\mathcal{M}.
 \end{equation}

\begin{obs}
In \cite{EGST1}, \cite{EGST2} the authors 
describe, using quiver Hopf algebras, the representation theory of the Drinfeld double $D(T_{n,m}^{*\cop})$ of the dual of the 
(co-opposite) generalized Taft algebra $T_{n,m}$.
Among other things, they prove that these are algebras of tame
representation type, present the indecomposable objects and describe the tensor product of two simple objects.
Since $_{D(T^{*\cop}_{n,m})}\mathcal{M} \simeq \ _{T_{n,m}^{*}}^{T_{n,m}^{*}}\mathcal{YD}\simeq\ _{T_{n,m}}^{T_{n,m}}\mathcal{YD} \simeq\
_{D}\mathcal{M}$, one possible way to describe $_{D}\mathcal{M}$ would be to
translate all the results using the cocycle deformation and the braided equivalences above. 
Another way is to use recent results by Pogorelsky and Vay to describe the simple modules and some projective envelopes 
using generalized Verma modules on weights.
Alternatively, we have
chosen to characterize directly the simple objects of $_{D}\mathcal{M}$, since we are particularly interested in the 
explicit structure of the modules, in order to compute their Nichols algebras.
\end{obs}

\subsection{Simple $D$-modules}\label{subsec:simpleDmodules}
In this subsection, we describe the simple left $D$-modules. 

For $0\leq i,j\leq nm-1$, let 
$r_{ij}$ be the positive integer satisfying that $1\leq r_{ij}\leq n$ and 
$$r_{ij}= \begin{cases}
               \begin{array}{lr}
                i+\frac{j}{m} +1 \mod n &\text{ if }\quad m \mid j,\\
                n &\text{ if }\quad m \nmid j.
               \end{array} 
               \end{cases}$$
We write simply $r=r_{ij}$ if no emphasis in $i,j$ is needed.

\begin{defi}\label{def:simplemodulesD} 
Let $\V$ be the  
$r$-dimensional $\ku$-vector space with basis $B=\{v_{0},\cdots, v_{r-1}\}$ endowed with a $D$-action determined by 
$$A\cdot v_{k} =  \xi^{i-k}v_{k}\qquad\qquad  g\cdot v_{k}= \xi^{j-km}v_{k} \qquad \forall\ 0\leq k\leq r-1,$$
$$x\cdot v_{k}= \begin{cases}
\begin{array}{lr}
 v_{k+1} & \text{ if }\quad 0\leq k<r-1,\\
 (1-\xi^{jn})v_{0} & \text{ if }\quad k=r-1,
\end{array} 
              \end{cases}
              $$
$$
              X\cdot v_{k}= \begin{cases}
\begin{array}{lr}
 0 & \text{ if }\quad k=0,\\
c_{k}v_{k-1} & \text{ if }\quad 0< k\leq r-1,\\
\end{array}              
\end{cases}
$$
where
\begin{align}\label{SolCk}
c_{k}=(k)_\omega\, \omega^{-k}(\xi^{j}\omega^{-k+1+i}-1), \qquad \forall\ 1\leq k\leq r-1.
\end{align}
\end{defi}

A direct but tedious calculation using the relations in Proposition \ref{prop:drinfeld_double} shows that 
the action of $D$ on $\V$ is well-defined. 
\begin{rmks}\label{rmk:simpleVij}{\rm
$(a)$ If $r<n$, then $(1-\xi^{jn})=0$, as $m\mid j$.  In particular, $x\cdot v_{r-1}=0$ and 
$x$ acts nilpotently on $\V$. Note that if $j=(n-i)m$, one has that $\dim \V = 1$.

$(b)$ Observe that $c_{k}\neq 0$ for all $1\leq k\leq r-1$. In fact, if $c_{k}=0$ then $\xi^{j + m(-k+1+i)}=1$, which implies that 
$j + m(-k+1+i) \equiv 0 \mod nm$. In such a case, $m\mid j$ and hence $k= \frac{j}{m} + i+1 = r$,  a contradiction.
%
%
%
} 
\end{rmks}

 \begin{lema}\label{lem:Vij-simple} The $D$-modules
 $\V$, $0\leq i,j< nm$, are simple and pairwise non-isomorphic. 
 \end{lema}

\pf Fix $0\leq i,j \leq nm-1$ and 
let $V$ be a $D$-submodule of $\V$. Take $0\neq v \in V$ and write $v=\sum_{\ell=0}^{r-1} \alpha_{\ell}v_{\ell}$ with $k=\max\{\ell:\ \alpha_{\ell}\neq 0\}$.
Then, $X^{k-1}\cdot v = \alpha_{k}c_{k}c_{k-1}\cdots c_{1}v_{0}$. Since by Remark \ref{rmk:simpleVij} $(b)$, one has that $c_{\ell}\neq 0$ for all 
$1\leq \ell \leq r-1$, it follows that $v_{0} \in V$. By acting with $x$ on $v_{0}$, one sees that $v_{\ell} \in V$
for all $0\leq \ell \leq r$, and consequently $V=\V$.

We show now that the family $\{\V\}_{1\leq i,j<nm}$ consists of pairwise non-isomorphic simple $D$-modules. 
Suppose  that $T:\V\rightarrow \mathcal{V}_{k,\ell}$ 
is an isomorphism of $D$-modules and denote by $B_{\V}=\{v_{t}\}_{0\leq t\leq r-1}$ and $B_{\mathcal{V}_{k,\ell}}=\{u_{t}\}_{0\leq t\leq s-1}$,
the linear bases of $\V$ and $\mathcal{V}_{k,\ell}$, respectively. Clearly, we must have $r=s$.  
Since  $T(v_{0})$ is an eigenvector for $g$, there exist $0\leq a\leq r-1$ and 
$\lambda\in \Bbbk^{\times}$ such that $T(v_{0}) =\lambda u_{a}= \lambda x^a u_{0}$. If $a>0$, then 
$0=T(X\cdot v_{0}) =\lambda X \cdot  u_{a}=\lambda c_{a} u_{a-1}$, 
where $c_{a}$ is the scalar given in \eqref{SolCk}. But this is impossible, since $c_{a}\neq 0$ by Remark \ref{rmk:simpleVij} $(b)$. 
Hence $a=0$ and consequently $T(v_{b}) = T(x^b v_{0})= \lambda\ x^b u_{0} = \lambda\ u_{b}$ for all $0\leq b\leq r-1$. This implies that $i=k$ and $j=\ell$.
\epf

The following theorem gives the classification of the simple $D$-modules. We prove it through several lemmata.
\begin{teor}\label{theo:RepSimples}
The family $\{\V\}_{0\leq i,j< nm}$ form a complete set of isomorphism classes of simple $D$-modules. 
In particular, there exist $(nm)^2$ non-isomorphic simple $D$-modules.  
\end{teor}

Let $V$ be a simple $D$-module and
consider
the linear
subspace $W= \ker X \subset V$. Since $X$ is a nilpotent element and
$\xi X A=AX$, $\omega  X U=UX$, $\omega  X g=gX$, it follows that $W$ is a non-zero subspace of $V$ 
invariant under the action of $g$, $A$ and $X$. 
Since $g$ commutes with $A$ and $U=A^{m}$, there exists a simultaneous eigenvector 
$v\in W$  of $A, U$ and $g$; say
\[g v=\alpha v, \quad A v=\gamma  v \quad\text{ and } \quad U v=\beta v, \,\,\, \text{with}\,\,\,\alpha,\gamma,\beta = \gamma^{m}\in \Bbbk.\]
Moreover, as $g^{nm}=A^{nm}=1$, we have that there exist $0\leq i,j\leq nm-1$ such that
$$
 A v=\xi^i  v\quad g v=\xi^j v,\quad U v=\omega^i v.
$$
%
%

\begin{lema} Let $V_{i,j}$ be the 
vector subspace spanned by $\{ v, xv, \ldots, x^{n-1} v \}$. Then
$V_{i,j}= V$. 
\end{lema}

\pf
To prove the claim, it is enough to show that $V_{i.j}$ is a $D$-submodule of V, \textit{i.e.} 
$V_{i,j}$
is invariant under the action of the generators $U$, $g$, $x$, $X$ and $A$. 
First we note that  $x^{k-1}v\in \ker X^{k}$, for all $0< k\leq n$.
Indeed, the result is true for $k=1$. Let $k>1$ and assume that $x^{k-2}v\in \ker X^{k-1}$. By \eqref{id-xX}, we have that 
\begin{align*}
 xX^k(x^{k-2}v)&=w^kX^k(x^{k-1}v)+(k)_\omega X^{k-1}\Gamma(k-1)(x^{k-2}v)\\
                 &=w^kX^k(x^{k-1}v)+(k)_\omega \Gamma(1-k)X^{k-1}(x^{k-2}v),
                \end{align*}
and consequently $X^k(x^{k-1}v)=0$. 

Since 
$U x^{k} v=\omega^{-k}x^kUv=\beta \omega^{-k} x^{k} v$ and $g x^{k} v=\omega^{-k}x^kgv=\alpha \omega^{-k} x^{k} v$
for all $0\leq k\leq n-1$, we have that 
$V_{i,j}$ is invariant under the action of $U$ and $g$. As $x^n=1-g^n$, it is clear that $V_{i,j}$ is also invariant under the action of $x$. Moreover, 
since by \eqref{id-xX1}, $Xx^{k}v$ belongs to the linear span of $ \{x^{k-1}v\}$, we have that $V_{i,j}$ is also $X$-invariant.
Finally, to show that $V_{i,j}$ is $A$-invariant, it is enough to prove that $A(x^{k}v)=\gamma \xi ^{-k}\,  x^{k}v$, for all $0\le k\leq n-1$.  
For this, we proceed by induction on $k$. The case $k=0$ follows from the equality $Av=\gamma  v$. 
Consider now $0< k< n$ 
and assume that $A(x^{k-1}v)=\gamma \xi ^{-(k-1)} \ x^{k-1}v$. 
Since  $X^{n-1}A=\xi^{-(n-1)}AX^{n-1}$ and $X^{n-1}(x^{k-1}v)=0$ for all $1\le k< n$,
it follows that $X^{n-1}(Ax^{k-1}v)=0$. Thus 
\begin{align*}
\xi A(x^{k}v)&=	\xi Ax(x^{k-1}v)\\
&=xA(x^{k-1}v)+\frac{1-\xi^n}{(n-1)!_{\omega}}
\left(gX^{n-1}AU-X^{n-1}A\right)(x^{k-1}v)\\
&=xA(x^{k-1}v)=x(\gamma \xi ^{-(k-1)}  x^{k-1}v)=\gamma \xi ^{-(k-1)} \ x^{k}v.
\end{align*}
Hence $A(x^{k}v)=\gamma \xi ^{-k} \ x^{k}v$, and the lemma is proved. 
\epf

Note that by \eqref{id-xX1}, for any $0\leq k\leq n-1$, there exists $d_{k}\in \Bbbk$ such that $X x^k v=d_{k}x^{k-1}v$. 
In fact, it turns out that $d_{k}= c_{k}$ for all $0\leq k\leq n-1$. Indeed, using again \eqref{id-xX1} we get
\begin{align*}
d_{k-1}x^{k-1} v&= x X (x^{k-1} v)=\omega Xx^{k} v+(1-Ug)x^{k-1}v\\&
=\omega\,  d_{k}x^{k-1} v+(1-\xi^{j} \omega^{i-2(k-1)}) x^{k-1} v\\
&
=(\omega\,  d_{k}+1-\xi^{j} \omega^{i-2(k-1)}) x^{k-1} v.
\end{align*}
Thus, 
$d_{k-1}=\omega\,  d_{k}+1-\xi^{j} \omega^{i-2(k-1)}$. 
On the other hand, we have that 
\begin{align*}
0&= x X ( v)=\omega X x v+(1-Ug) v
=\omega\,  d_{1} v+ (1-\xi^{j} \omega^{i})v=(\omega\,  d_{1}+1-\xi^{j} \omega^{i})v,
\end{align*}
which implies that $d_{1}=\omega^{-1}(\xi^{j} \omega^{i}-1)=c_{1}$. Consequently   
$$d_{k}= (k)_{\omega^{-1}}(\omega^{-k}\xi^{j} \omega^{i}-\omega^{-1})\,=(k)_\omega\, \omega^{-k}(\omega^{-(k-1)}\xi^{j} \omega^{i}-1) = c_{k}
\quad \text{ for all }0\leq k\leq n-1.$$

\begin{lema}\label{lem:dimension}
For all $0\leq i,j\leq nm-1$, it holds that 
$\dim V_{i,j}= r_{ij}$.
\end{lema}

\pf
As $U v =\omega^i v$, we have that  $U x^k v=\omega^{i-k} x^k v $ 
for all $0\leq k\leq n-1$. Set $r=\max\{1\leq k\leq n:\ x^{k-1} v \neq 0\}$.
Then,
$v, x v, \ldots x^{r-1} v $ are eigenvectors of $U$ associated with distinct eigenvalues, and consequently
the set $\{ v, x, \ldots, x^{r-1} v \}$ is a basis of $V_{i,j}$.

If $m\nmid j$, then $1\neq \xi^{jn}$ and whence $x^{n}$ acts non-trivially on $V_{i,j}$. This implies that 
$\{v,xv,\ldots, x^{n-1}v\}$ is a basis of $V_{i,j}$ and $\dim V_{i,j}= n = r_{ij}$.

Assume now that $m\mid j$. Then 
$x^{n}\cdot v = (1-g^{n})\cdot v = (1-\xi^{jn})v=0$ and $x^{n}$ acts by 0 on $V_{i,j}$. 
Since $x^r$ also acts by zero on $V$, we have that  
$$
c_{r-1}x^{r-1} v = x X (x^{r-1} v)=\omega Xx^{r} v+ (1-Ug)x^{r-1} v
= (1-\xi^j\omega^{i}\omega^{-2(r-1)}) x^{r-1} v.
$$
Comparing with  \eqref{SolCk} we have that
$$
\omega^{-(r-1)}\xi^j\omega^{i}\left(\omega^{-(r-1)}+(r-1)_\omega\, \omega^{-(r-2)}\right) =
1+(r-1)_\omega\, \omega^{-(r-1)}.
$$
But
$
1+(r-1)_\omega\, \omega^{-(r-1)}$ $=1+(1+\omega+\cdots+\omega^{r-2})\omega^{-(r-1)}$ 
and
\begin{align*}
\left(\omega^{-(r-1)}+(r-1)_\omega\, \omega^{-(r-2)}\right)&=\omega^{-(r-1)}+(1+\omega+\cdots+\omega^{r-2})
\omega^{-(r-2)}\\&=1+\omega^{-1}+\cdots+\omega^{-(r-1)}=(r)_{\omega^{-1}},
\end{align*}
which implies that $\omega^{-(r-1)}\xi^j\omega^{i}\  (r)_{\omega^{-1}} = (r)_{\omega^{-1}}$.
If $r<n$, then $ (r)_{\omega^{-1}}\neq 0$ and whence $\omega^{r-1}= \xi^j\omega^{i}=\omega^{i+\frac{j}{m}}$. 
Hence $r= i +\frac{j}{m} + 1 \mod n=r_{i,j}$. 
Suppose now that $r=n$. Then 
\begin{align*}
(\xi A x-x A ) x^{n-1} v&=\xi A x^{n} v- \xi^i\xi^{-(n-1)} x^{n}v
=\xi A(1-g^n)v-\xi^i\xi^{-(n-1)}(1-g^n)v\\
&=(1-\xi^{jn})\xi^i(\xi-\xi^{-(n-1)})v=0.
\end{align*}
Observe also that
\begin{equation}\label{c1c2..cn-1}
\begin{aligned}
X^{n-1} x^{n-1} v&=X^{n-2}(X x^{n-1} v)=c_{n-1}X^{n-2} x^{n-2} v
                   =c_{n-1}c_{n-2}X^{n-3} x^{n-3} v\\
                   &\hspace{2cm}\vdots \\&=c_{n-1}c_{n-2}\cdots c_{1} v
                   =(n-1)!_{\omega}\prod_{k=1}^{n-1}\omega^{-k}(\omega^{-(k-1)}\xi^j\omega^i-1)v\\
                       &=(n-1)!_{\omega}\ \omega^{-\frac{(n-1)n}{2}}\prod_{k=0}^{n-2}(\omega^{-k}\xi^j\omega^i-1)v\\
                       &=(n-1)!_{\omega}\omega^{-\frac{(n-1)n}{2}}(-1)^{n-1}\prod_{k=0}^{n-2}(1-\omega^{-k}\xi^j\omega^i)v\\
                       &=(n-1)!_{\omega}(-1)^{n-1}(-1)^{n-1}\prod_{k=0}^{n-2}(1-\omega^{-k}\xi^j\omega^i)v\\
                       &=(n-1)!_{\omega}\prod_{k=0}^{n-2}(1-\omega^{-k}\xi^j\omega^i)v.
\end{aligned}
\end{equation}
Writing $\lambda_j=\displaystyle\prod_{k=0}^{j}(1-\omega^{-k}\xi^j\omega^i)\in \Bbbk$, we have that
\begin{align*}
(gX^{n-1}AU-X^{n-1}A)(x^{n-1} v)&=gX^{n-1}AU(x^{n-1} v)-X^{n-1}A (x^{n-1} v)\\
& =\xi^i \xi^{-(n-1)}(\omega^i \omega^{-(n-1)}g-1)(X^{n-1} x^{n-1} v)\\
&=\xi^i \xi^{-(n-1)}(\omega^i \omega^{-(n-1)}g-1)((n-1)!_{\omega}\lambda_{n-2}v) \\
&=\xi^i \xi^{-(n-1)}(n-1)!_{\omega}\lambda_{n-2}(\xi^j\omega^i \omega^{-(n-1)}-1)v\\
&=-\xi^i \xi^{-(n-1)}(n-1)!_{\omega}\lambda_{n-1}v.
\end{align*}
Thus, \eqref{eqdouble3} implies that
\begin{align*} (1-\xi^{jn})\xi^i(\xi-\xi^{-(n-1)})v&=\frac{1-\xi^n}{(n-1)!_{\omega}}\left(-\xi^i \xi^{-(n-1)}(n-1)!_{\omega}\lambda_{n-1}\right)v\\
&=\xi^i(\xi-\xi^{-(n-1)})\lambda_{n-1}v.
\end{align*}
Then
$
0=(1-\xi^{jn})=\prod_{k=0}^{n-1}(1-\omega^{-k}\xi^j \omega^i)=\lambda_{n-1}
$.
Hence, there exists $k$ with $0\le k\le n-1$ such that 
$1-\omega^{-k}\xi^j \omega^i=0$. Since by Remark \ref{rmk:simpleVij} $(b)$, $c_{k}\neq 0$ for all $1\leq k\leq n-1$,
it follows from \eqref{c1c2..cn-1} that $1-\omega^{-k}\xi^j \omega^i\neq 0$ for all $0\leq k\leq n-2$. Thus $1-\omega^{-(n-1)}\xi^j \omega^i=0$ 
and consequently $\omega^{-1}=\xi^j \omega^i=\omega^{\frac{j}{m}+i}$, which implies that $n= i+\frac{j}{m} +1 \mod n = r_{ij}$. 
\epf

\begin{lema}\label{lem:Vij-V}
For all $0\leq i,j\leq nm-1$, it holds that 
$V_{i,j}\simeq  \V$ as $D$-modules. 
\end{lema}

\pf Define the linear map $\varphi: V_{i,j} \to \V$ by $\varphi(xv^{k}) = v_{k}$ for all $1\leq k\leq r-1$. Using the calculations
above, it is easy to see that $\varphi$ is indeed an isomorphim of $D$-modules.
\epf

 \noindent \textit{Proof of Theorem \ref{theo:RepSimples}.} 
 Follows directly from Lemmata \ref{lem:Vij-simple}, \ref{lem:dimension} and \ref{lem:Vij-V}.

\begin{defi}\label{def:grafo-modulo}
Let $V$ be a $D$-module with linear basis $B=\{v_1,\ldots,v_t\}$.  
The oriented graph associated with $V$, with respect to the basis $B$, is the graph whose
vertices are the vectors $v_{j}$ with $1\leq j\leq t$, and the arrows are given as follows: we draw an arrow from $v_{\ell}$ to $v_{k}$ 
if $x v_\ell=\sum_{i=1}^t\lambda_{i,\ell} v_i$ with $\lambda_{k,\ell}\neq 0$, and we draw an arrow from $v_{\ell}$ to $v_{k}$ 
if $X v_\ell=\sum_{i=1}^t\lambda^{\prime}_{i,\ell} v_i$ with $\lambda^{\prime}_{k,\ell}\neq 0$. 
\end{defi}

\begin{example}\label{RemarkGraph} {\rm Let $\V$ be a simple $D$-module. 
When $x$ acts as a nilpotent element on $\V$, the oriented graph 
representing its action is the following
\[{\tiny\xymatrix@C=5mm{v_0  \ar@{->}@/^/[r] 
		& v_1  \ar@{->}@/^/[r]\ar@{->}@/^/[l]
		& v_2\ar@{->}@/^/[r]\ar@{->}@/^/[l]&\ldots \ar@{->}@/^/[r] \ar@{->}@/^/[l] & 
		v_{r-3}\ar@{->}@/^/[r] \ar@{->}@/^/[l] & v_{r-2}\ar@{->}@/^/[r]\ar@{->}@/^/[l]	
		& v_{r-1}  \ar@{->}@/^/[l] }}\]
If $x^{n}$ acts as a semisimple element, the oriented graph is  

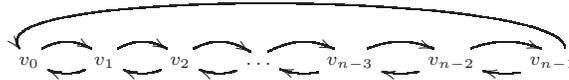
\bigbreak

\[{\tiny\xymatrix@C=5mm{v_0  \ar@{->}@/^/[r] 
		& v_1  \ar@{->}@/^/[r]\ar@{->}@/^/[l]
		& v_2\ar@{->}@/^/[r]\ar@{->}@/^/[l]&\ldots \ar@{->}@/^/[r] \ar@{->}@/^/[l] & 
		v_{n-3}\ar@{->}@/^/[r] \ar@{->}@/^/[l] & v_{n-2}\ar@{->}@/^/[r]\ar@{->}@/^/[l]	
		& v_{n-1}  \ar@{->}@/^/[l] \ \ \ar@{->}@(ur,ul)[llllll]!<-1mm,1mm>}}\]
		}
\end{example}

\bigbreak

\subsection{Projective covers of simple objects}
In this subsection we compute the projective covers of the simple $D$-modules $\V$ 
given in \S\S \ref{subsec:simpleDmodules}. We write 
$ \mathcal{P}(\V)$ for the projective cover of $\mathcal{V}_{i,j}$.
We study the cases $\dim \V=n$ and $\dim \V<n$ separately.

We start by observing that, for all $0\leq i,j\leq nm-1$, the elements 
$$
e_{i,j}=\frac{1}{nm}\sum_{r,s=0}^{nm-1} \xi^{-ir-js}A^rg^s
$$
form a complete set of non-zero orthogonal idempotents of $D$, that is $\sum_{i,j=0}^{nm-1} e_{i,j}=1$ and 
$e_{i,j}e_{k,\ell} = \delta_{i,k}\delta_{j,\ell}$ for all $0\leq i,j,k,\ell\leq nm-1$. Moreover, one has that 
\begin{align}\label{comute-idem}
&ge_{i,j}=\xi^je_{i,j},& Ae_{i,j}=\xi^ie_{i,j},& & \Gamma(t)e_{i,j}=\gamma(t)e_{i,j},&
\end{align}
where $\gamma(t):=1-\xi^{m(t+i)+j}$ and $\Gamma(t)$ is given as in \eqref{ref:lema29}, for any integer $t$. 
Since $X^kx^kg=gX^kx^k$, it follows from \eqref{id-xkXkA} that 
\begin{align}\label{eq:xkXkA}
X^kx^ke_{i,j}=e_{i,j}X^kx^k,\qquad\text{for all}\,\,\, k\geq 1.&
\end{align}

\begin{teor}\label{teor:cover-n}
Let  $\mathcal{V}_{i,j}$ be a simple $D$-module with $\dim \V = n$. Then
$\mathcal{P}(\V)\simeq \V$.	
\end{teor}

\pf We prove that $\V$ is a projective $D$-module by showing that it is isomorphic as $D$-module to $De^{(n)}_{i,j}$,
with $e^{(n)}_{i,j}=\lambda_n X^{n-1}x^{n-1}e_{i,j}$ an idempotent element of $D$. 

Set $\lambda_n:=1/(c_{1}\ldots c_{n-1})$, with $c_{k}=-\gamma(1-t)w^{-t}(t)_{\omega}$ the 
coefficients given in 
\eqref{SolCk}. Note that,
by Remark \ref{rmk:simpleVij} $(b)$, these are non-zero. 
Applying \eqref{id-xX1} recursively, one can see
that $X^{n-1}x^{n-1}X^{n-1}=(n-1)!_{\omega}X^{n-1}\Gamma(0)\Gamma(1)\ldots \Gamma(n-2)$. Using this, we have that
\begin{align*}
e^{(n)}_{i,j}e^{(n)}_{i,j}&\overset{\eqref{eq:xkXkA}}{=}\lambda_n^2X^{n-1}\left(x^{n-1}X^{n-1}\right)x^{n-1}e_{i,j}\\
                          &=\lambda_n^2(n-1)!_{\omega}X^{n-1}\Gamma(0)\Gamma(1)\ldots \Gamma(n-2)x^{n-1}e_{i,j}\\
                          &\overset{\eqref{id-gammax}}{=}\lambda_n^2(n-1)!_{\omega}X^{n-1}x^{n-1}\Gamma(2)\Gamma(3)\ldots \Gamma(n-1)\Gamma(0)e_{i,j}\\
                          &\overset{\eqref{comute-idem}}{=} \lambda_n^2(n-1)!_{\omega}\gamma(0)\gamma(2)\gamma(3)\ldots \gamma(n-1)X^{n-1}x^{n-1}e_{i,j}\\
                          &= \lambda_nX^{n-1}x^{n-1}e_{i,j}
                          =e^{(n)}_{i,j}.
\end{align*}
In order to prove that $\V\simeq De^{(n)}_{i,j}$, we show that the set
$B=\{e^{(n)}_{i,j},xe^{(n)}_{i,j},\ldots,x^{n-1}e^{(n)}_{i,j}\}$ is a linear basis of $De^{(n)}_{i,j}$. Then, the  
linear map $\psi:De^{(n)}_{i,j}\to \mathcal{V}_{i,j}$ given by $\psi(x^ke^{(n)}_{i,j})=v_{k}$ for all $0\leq k\leq n-1$, yields the desired
isomorphism.
Note first 
that $Xe^{(n)}_{i,j}=0$ and $gX^{n-1}x^{n-1}=X^{n-1}x^{n-1}g$. Thus, by \eqref{comute-idem} we get that  $ge^{(n)}_{i,j}=\xi^je^{(n)}_{i,j}$. 
Analogously, from \eqref{id-xkXkA} and \eqref{comute-idem} we obtain $Ae^{(n)}_{i,j}=\xi^i e^{(n)}_{i,j}$. Moreover,
for all $1\leq k\leq n-1$ we have
\begin{align*}
X x^ke^{(n)}_{i,j}&\overset{\eqref{id-xX1}}{=}\omega^{-k}\left(x^kX-(k)_{\omega}\Gamma(k-1)x^{k-1}\right) e^{(n)}_{i,j}
=-\omega^{-k}(k)_{\omega}\Gamma(k-1)x^{k-1} e^{(n)}_{i,j}\\
                  &\overset{\eqref{id-gammax}}{=} -\omega^{-k}(k)_{\omega}x^{k-1} \Gamma(k-1-2(k-1)) e^{(n)}_{i,j} \\
                  & \overset{\eqref{comute-idem}}{=}-\omega^{-k}(k)_{\omega}\gamma(1-k)x^{k-1}e^{(n)}_{i,j}=c_{k}\ x^{k-1}e^{(n)}_{i,j}.
\end{align*}
This implies that $B$ is a linearly independent set and generates $De^{(n)}_{i,j}$.
\epf

Next we proceed to compute the projective covers of the simple modules 
$\V$ with $\dim \V <n$. Let $1\leq i,j <nm$ be such that $m\mid j$ and $r= i+\frac{j}{m} +1 \not\equiv 0\mod n$.
Consider the $2n$-dimensional $\Bbbk$-vector space $M_{i,j}$ with linear basis  
\[B=\{u^{(1)}_0,\ldots,u^{(1)}_{n-r-1},\, v^{(1)}_0, \ldots, v^{(1)}_{r-1},\,u^{(2)}_0,\ldots,u^{(2)}_{n-r-1},\,v^{(2)}_0,\ldots,v^{(2)}_{r-1}  \}.\] 
We define a $D$-module structure on $M_{i,j}$ as follows.  The actions of $g$ and $A$ are given by
\begin{equation*}
\left\{
\begin{aligned}
& gu^{(t)}_{k}=\omega^{j-r-k}u^{(t)}_{k},\qquad gv^{(t)}_{\ell}=\omega^{j-\ell}v^{(t)}_{\ell},\qquad t=1,2,  \\[2pt]
&Au^{(1)}_{k}=\xi^{n+i-r-k} u^{(1)}_{k},\quad Au^{(2)}_{k}=\xi^{i-r-k}(u^{(2)}_{k}+u^{(1)}_{k}),\\[2pt]
&Av^{(t)}_{\ell}=\xi^{i-\ell} v^{(t)}_{\ell},\qquad\qquad t=1,2,
\end{aligned}
\right.
\end{equation*}
for all $0\leq k\leq n-r-1$ and  $0\leq \ell\leq r-1$; whereas
the action of $x$ is given by
\begin{equation*}
\left\{
\begin{aligned}
&xu^{(t)}_{k}=u^{(t)}_{k+1},\quad\,\, xu^{(1)}_{n-r-1}=v^{(1)}_0,\quad xu^{(2)}_{n-r-1}=0,\quad t=1,2,  \\[2pt]
&xv^{(t)}_{\ell}=v^{(t)}_{\ell+1},\qquad xv^{(2)}_{r-1}=u^{(2)}_0,\qquad xv^{(1)}_{r-1}=0,  \,\,\,\qquad t=1,2,
\end{aligned}
\right.
\end{equation*}
for all $0\leq k\leq n-r-2$ and $0\leq \ell\leq r-2$; and the action of $X$ is defined by
\begin{equation*}
\left\{
\begin{aligned}
&Xu^{(t)}_{k}=c_{k}u^{(t)}_{k-1},\quad Xu^{(2)}_{0}=\alpha\omega^{-r} v^{(1)}_{r-1},\qquad \quad Xu^{(1)}_{0}=Xv^{(1)}_{0}=0,  \\[2pt]
&Xv^{(1)}_{\ell}=c_{\ell}v^{(1)}_{\ell-1}, \quad\,\,\,\, Xv^{(2)}_{\ell}=\alpha\omega^{-\ell} v^{(1)}_{\ell-1}+c_{\ell}v^{(2)}_{\ell-1},\quad  
Xv^{(2)}_{0}=\alpha u^{(1)}_{n-r-1},
\end{aligned}
\right.
\end{equation*}
where $\alpha\in \ku^{\times}$ satisfies that
$$\alpha(\omega^{-r}-1)\frac{1-\xi^n}{(n-1)!_{\omega}}\prod_{a=0}^{r-1}c_a\prod_{b=0}^{n-r-1}c_{b}=1,$$
$t=1,2$, $\,1\leq k\leq n-r-1$ and $1\leq \ell\leq r-1$. 
Note that $\ker X|_{M_{i,j}}=\Bbbk\{ u^{(1)}_0,v^{(1)}_0\}$.

The oriented graph of $M_{i,j}$ associated with the basis $B$ is as follows

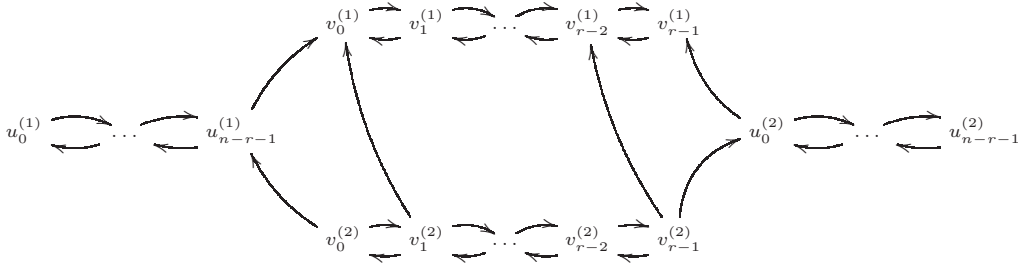
\begin{tiny}
$$\xymatrix@C=4mm{
		&   &   & v^{(1)}_0\ar@{->}@/^/[r]
		& v^{(1)}_1 \ar@{->}@/^/[r]\ar@{->}@/^/[l]
		& \ldots\ar@{->}@/^/[r]\ar@{->}@/^/[l]
		& v^{(1)}_{r-2} \ar@{->}@/^/[r]\ar@{->}@/^/[l]
		& v^{(1)}_{r-1}\ar@{->}@/^/[l]
		&\\		
		u^{(1)}_0 \ar@{->}@/^/[r] 
		& \quad \ldots\quad \ar@{->}@/^/[r]\ar@{->}@/^/[l]
		& u^{(1)}_{n-r-1}\ar@{->}@/^/[ur]\ar@{->}@/^/[l] 
		&
		& 
		& 
		&
		& 
		& u^{(2)}_0  \ar@{->}@/^/[r]\ar@{->}@/^/[ul]+<1mm,-2mm>
		&  \quad \ldots\quad  \ar@{->}@/^/[r]\ar@{->}@/^/[l] 
		& u^{(2)}_{n-r-1}\ar@{->}@/^/[l]\\ 
		&   &   & v^{(2)}_0\ar@{->}@/^/[r]\ar@{->}@/^/[ul]
		& v^{(2)}_1  \ar@{->}@/^/[r]\ar@{->}@/^/[l]\ar@{->}@/^/[luu]
		& \ldots\ar@{->}@/^/[r]\ar@{->}@/^/[l]
		& v^{(2)}_{r-2} \ar@{->}@/^/[r]\ar@{->}@/^/[l]
		& v^{(2)}_{r-1}\ar@{->}@/^/[l]\ar@{->}@/^/[luu]\ar@{->}@/^1pc/[ur] 
		&}$$
 \end{tiny}

		
It is clear from the graph above that  
\[M_0=  \Bbbk\left\{ v^{(1)}_0, \ldots, v^{(1)}_{r-1}\right\}, \qquad M_1= \Bbbk\left\{ u^{(1)}_0,\ldots,u^{(1)}_{n-r-1},\, v^{(1)}_0, \ldots, v^{(1)}_{r-1} \right\}, \]
\[M_2=\Bbbk\left\{ u^{(1)}_0,\ldots,u^{(1)}_{n-r-1},\, v^{(1)}_0, \ldots, v^{(1)}_{r-1},\,u^{(2)}_0,\ldots,u^{(2)}_{n-r-1}\right\}, 
\quad\,\,\, M_3=M_{i,j}, \]
are $D$-submodules  of $M_{i,j}$. Moreover, these modules give us a composition series of $M_{i,j}$.

\begin{prop}\label{prop-serie}
	The composition series of $M_{i,j}$ is
$M_0 \subset M_1 \subset M_2 \subset M_3$,
with composition factors
\[M_0\simeq \mathcal{V}_{i,j} ,\quad M_1/M_0\simeq \mathcal{V}_{n+i-r,j-mr}, 
\quad M_2/M_1\simeq \mathcal{V}_{i-r,j-mr},\quad M_3/M_2\simeq \mathcal{V}_{i,j}.\] 
\end{prop}

\pf Follows by a direct computation.
\epf

\begin{prop} \label{prop:socle}
${\rm soc}(M_{i,j})=\V$. In particular, $M_{i,j}$ is an indecomposable $D$-module.
\end{prop}

\pf  Let $V$ be a simple $D$-submodule of $M_{i,j}$. By the 
proof of Theorem \ref{theo:RepSimples},  
there exists $0\neq v\in \ker X|_{V}$, which is a simultaneous eigenvector of $A$ and $g$, and satisfies that $V=Dv$. 
Thus, there are $\beta_1,\beta_2\in \Bbbk$ such that $v=\beta_1u^{(1)}_0+\beta_2v^{(1)}_0$. 
As $Av=\beta_1\xi^{n+i-r}u^{(1)}_0+\beta_2\xi^{i}v^{(1)}_0$ and 
$\xi^{n+i-r}\neq \xi^{i}$, we must have that 
$\beta_1=0$ or $\beta_2=0$. If $\beta_2=0$ then $V=Dv=M_1\supset \V$,
which implies that $V$ is not simple, a contradiction. 
Hence, $v=\beta_2v^{(1)}_0$ and $V=\V$.  
\epf

Denote by $\mathcal{I}(\V)$ the injective envelope of  $\V$. Note that, as
$D$ is a self-injective algebra (because it is a symmetric algebra), one has that $\mathcal{I}(\V)\simeq \mathcal{P}(\V)$ as $D$-modules. 
Also injective $D$-modules are projective and vice-versa. 

\begin{teor}\label{teo:cover} 
Let $\V$ be a simple $D$-module with $\dim \V <n$. Then $ \mathcal{P}(\V)\simeq M_{i,j}$.
\end{teor}

\pf 

Let $\iota:\V \to \mathcal{I}(\V)$ be the monomorphism given by the inclusion of $\V$ in its 
injective envelope and denote by $\kappa:\V\to M_{i,j}$ the inclusion given by Proposition \ref{prop-serie}.
Since $\mathcal{I}(\V)$ is injective, there exists a $D$-module map $\nu:M_{i,j} \to \mathcal{I}(\V)$ 
such that $\nu\kappa=\iota$. As $\V\simeq \kappa(\V)$ is a simple $D$-module, 
we have that $\ker \nu \cap \kappa(\V)=0$. Hence $\ker \nu=0$ and whence $\nu$ is an injective map. In particular, 
$2n=\dim M_{i,j} \leq \dim \mathcal{I}(\V)$. Since 
\[_{D}D\simeq 
\bigoplus_{0\le i,j\le nm-1}\,\,\mathcal{I}(\V)^{r}\\
=
\bigoplus_{r\neq n}\,\, \mathcal{I}(\V)^{r} \, \oplus \ 
\bigoplus_{r=n} \mathcal{I}(\V)^n, \]
as $D$-modules, and $\V\simeq \mathcal{P}(\V)\simeq \mathcal{I}(\V)$ for $\dim \V =n$, 
by Theorem \ref{teor:cover-n} we get that
\begin{align*}
\dim D=n^4m^2&\geq \sum_{r=1}^{n-1}nm\cdot r \cdot 2n+ nm(nm-n+1)\cdot n \cdot n \\
&=n^4m^2- n^4m+n^3m +n^2m n(n-1)= n^4m^2.
\end{align*}
This implies that necessarily $\dim \mathcal{I}(\V)=2n$, and $\nu$ is an isomorphism of $D$-modules.
\epf

\subsection{Braided equivalences}\label{sub-braided} 
In this subsection, on the goal to study Nichols algebras over the simple modules $\{\V\}_{0\leq i,j< nm}$, we 
describe the explicit relation between the simple $D$-modules and the simple objects in $ _{T_{n.m}}^{T_{n,m}}\mathcal{YD}$.

Denote by $F:\  _{D}\M \to\  _{T_{n.m}}^{T_{n,m}}\mathcal{YD}$ the braided tensor equivalence given by the 
composition of the functors given in Proposition \ref{prop:equiv-D-mod-YD}, Equation \eqref{eq:braided-equiv} and the beginning 
of Section \ref{sec:D-modules}:
$$\xymatrix{F:\ _{D}\M \ar[r]^{F_{1}}& \ydhnm \ar[r]^{F_{2}}& \ _{R_{n.m}}^{R_{n,m}}\mathcal{YD} \ar[r]^{F_{3}} 
& _{T_{n.m}}^{T_{n,m}}\mathcal{YD}.}
$$
Here
$F(V) = V$ as $\Bbbk$-vector spaces for any $V\in\, _{D}\M $.


\begin{prop}\label{prop:double-tnm} 
	For $0\leq i,j<  nm$, the object $F(\V) \in \ _{T_{n,m}}^{T_{n,m}}\mathcal{YD}$ is simple and its structure is determined by the 
	following equalities:
	\begin{align}\label{double-tnm1}
	&g\cdot v_\ell=\xi^{-j}\omega^\ell v_{\ell},\qquad\quad 
	 x\cdot v_\ell =
	\begin{cases}
\begin{array}{cl}
 - \xi^{-j}\omega^\ell v_{\ell+1}& \text{ if } \quad 0\leq \ell<r-1;\\
 0& \text{ if } \quad \ell=r-1;
\end{array}
	\end{cases}\\
\label{double-tnm2}  
	&\rho (v_\ell)=\sum_{k=0}^{\ell}\beta^{i,j}_{k,\ell}\,\,   x^k g^{-i+\ell-k} \otimes v_{\ell-k},
	\end{align}
for all $0\leq \ell \leq r-1$, where
$\beta^{i,j}_{k, \ell}={\ell \choose k}_{\omega} \prod_{s= \ell- k}^{\ell-1} \left(\omega^{-i}-\xi^j\omega^{-s}\right)$ if $k\geq 1$, and $\beta^{i,j}_{0,l}=1$.
\end{prop}

\pf The first assertion is obvious.  Let us check the Yetter-Drinfeld module structure of $F(\V)$ over $T_{n,m}$. To this end, we describe the 
images of $\V$ under the composition of the functors $F_{1}$, $F_{2}$ and $F_{3}$.
By Proposition \ref{prop:equiv-D-mod-YD}, we have that $F_{1}(\V) \in \ydhnm$ via
$$A\cdot v_{\ell} =  \xi^{i-\ell}v_{\ell}, \qquad
              X\cdot v_{\ell}= \begin{cases}
\begin{array}{cl}
 0 & \text{ if }\quad \ell=0,\\
c_{\ell}v_{\ell-1} & \text{ if }\quad 0< \ell\leq r-1,\\
\end{array}              
\end{cases}
$$
$$ \rho(v_{\ell}) = \sum_{k=0}^{r-1-\ell}\dfrac{1}{(k)_\omega!} X^k  A^j U^{-(k+\ell)} \otimes v_{k+\ell}+
\sum_{k=r-\ell}^{r-1}\dfrac{1-\xi^{jn}}{(k)_\omega!} X^k A^j U^{-(k+\ell)} \otimes v_{k+\ell-r},$$
for all $0\leq \ell \leq r-1$.
Besides, from \eqref{pairing-ACE} it follows that the linear basis $\{X^aA^b\}_{ 0\le a,b\le nm-1}$ of $H=R_{n,m}^*$ 
has dual basis $\{u_{a,b}\}_{ 0\le a,b\le nm-1}$ given by 
\[	u_{a,b}=\frac{x^a}{nm\,  (a)_\omega!}\sum_{d=0}^{nm-1}\xi^{-bd}g^d. \]
Using these bases, from \eqref{eq:dualyd} we get that the action of $R_{n,m}$ on 
$F_{2}(F_{1}(\V)) \in \ _{R_{n,m}}^{R_{n,m}}\mathcal{YD}$ is   
$$
 g \cdot v_\ell= \xi^{-j}\omega^\ell v_{\ell}, \qquad 
 x \cdot v_\ell= \begin{cases}
               \begin{array}{ll}
                 - \xi^{-j}\omega^\ell v_{\ell+1} &\text{ if }\quad 0< \ell\leq r-1,\\
                 - \xi^{-j}\omega^{-1}(1-\xi^{jn})v_{0}& \text{ if }\quad \ell=r-1;
               \end{array}
              \end{cases}$$
whereas the coaction is given by
\begin{align*}
\rho(v_\ell)=\sum_{k=0}^{\ell}\beta^{i,j}_{k,\ell}\,\,   x^k g^{-i+\ell-k} \otimes v_{\ell-k},\ 
\beta^{i,j}_{k,\ell}=
\begin{cases}
1,& \text{if } k=0,\\
\frac{(-1)^k}{(k)!} w^{-\frac{k(k-1)}{2}} w^{-k(i-\ell)} c_\ell\ldots c_{\ell-k+1},& \text{if } k>0,
\end{cases}
\end{align*}
for all $0\leq \ell \leq r-1$. Using induction, one may check that the coefficients  $\beta^{i,j}_{k,\ell}$ satisfy that
\begin{align*}
\beta^{i,j}_{k,\ell}={\ell \choose k}_{\omega}\,\, \prod_{s= \ell- k}^{\ell-1} \left(\omega^{-i}-\xi^j\omega^{-s}\right), \quad\text{ for } k>0.
\end{align*}

Finally, by \S \ref{subsec:Radford-alg} we know that the Radford algebra $R_{n,m}$ is isomorphic to a $2$-cocycle
deformation of the Taft algebra $T_{n,m}$, where the $2$-cocycle $\sigma=e^{\eta}=\eps\ot \eps + \eta$ 
is the one defined in \eqref{eq:formula-cocycle}. 
Then, we may describe the structure of $F(\V) \in \ _{T_{n,m}}^{T_{n,m}}\mathcal{YD}$ 
by using that the functor $F_{3}$ does not change the 
comodule structure, but  the action is changed by the $2$-cocycle $\sigma^{-1}=e^{-\eta}=\eps\ot \eps - \eta$ as follows:
$$
g \cdot_{\sigma^{-1}} v_\ell= \xi^{-j}\omega^\ell v_{\ell},\qquad 
x \cdot_{\sigma^{-1}} v_\ell =
\begin{cases}
\begin{array}{lr}
- \xi^{-j}\omega^\ell v_{\ell+1} &\text{ if }\quad 0\leq \ell < r-1,\\
0  &\text{ if }\quad \ell=r-1,
\end{array}
\end{cases}
$$
for all $0\leq \ell \leq r-1$. This proves the proposition.
\epf


\section{On finite-dimensional Nichols algebras associated with $\V$ }\label{sec:Nicholsalg}
The aim of this section is to describe families of finite-dimensional Nichols algebras over some simple
objects $\V$ of the category of Yetter-Drinfeld modules over $H_{n,m}$, in the case when $n=2$ and $m\geq 2$, see Theorem \ref{prop:classifica}.
To this end,  
we first establish
the pairs $(i,j)$ such that 
$\dim \toba(\V)$ is finite by using recent results of Andruskiewitsch and Angiono \cite{AA}, and the 
classification of Heckenberger \cite{H2} of arithmetic root systems of rank 2. Then we describe 
by generators and relations certain families of Nichols algebras in small rank, see \S \ref{subsec:presentationNichols}.
We keep the notation used in the previous sections.

\subsection{Simple Yetter-Drinfeld modules over basic Hopf algebras}
We begin this section by presenting another way to describe the simple objects of $\, _{T_{n,m}}^{T_{n,m}}\mathcal{YD}$;
it is based on \cite[$\S\,3.3$]{AHS} and \cite{AA}. 

Fix $g$ and $\chi$ generators of the cyclic groups $C_{nm}$ and $\widehat{C_{nm}}\simeq C_{nm}$, respectively, such that 
$\chi(g)=\xi$.
Recall that $T_{n,m}=\toba(V) \# \ku C_{nm}$, where  $V=\ku\{x\} \in\  _{\ku C_{nm}}^{\ku C_{nm}}\mathcal{YD}$ 
is such that by
$g\cdot x=\omega x$, $\rho(x)=g\otimes x$, and 
$\toba(V) = \ku[x]/(x^{n})$ is the truncated polynomial algebra.

The category
$\, _{\ku C_{nm}}^{\ku C_{nm}}\mathcal{YD}$ is  semisimple and its simple objects are parametrized by pairs of natural numbers
$(a,b)$ with $0\leq a,b\leq nm-1$. Precisely,
the simple object $\lambda_{a,b}$ associated with $(a,b)$ 
is the one-dimensional vector space $V_{g^{a}}^{\chi^{b}}= \ku \{y\}$ generated by the element $y$ with action and coation given respectively by:
\begin{align}\label{struc-lambdaab}
g\cdot y=\xi^b y,&&\rho(y)=g^a\otimes y.
\end{align}

Consider now  $W_{a,b}:=V\oplus \lambda_{a,b} \in \, _{\ku C_{nm}}^{\ku C_{nm}}\mathcal{YD}$. 
By \cite[$\S\,3.3$]{AHS}, we have that its Nichols algebra $\toba(W_{a,b})$ is isomorphic to a braided bosonization 
$\toba(W_{a,b})\simeq \mathcal{K}\# \toba (V)$, where  $\mathcal{K}$ is a Hopf algebra in 
$_{T_{n,m}}^{T_{n,m}}\mathcal{YD}$. Here, the action of $T_{n,m}$ on $\mathcal{K}$ is given by the induced action of $g$ by \eqref{struc-lambdaab}
and the braided
adjoint action of $x$. In our case, the latter reads
$x\cdot z=\operatorname{ad}_c x(z)= xz - (g\cdot z) x$ for all $z\in \mathcal{K}$.
In particular, we have that $x\cdot y= xy - \xi^{b}y x$.

For all $0\leq a,b\leq nm-1$, lets write $L(\lambda_{a,b})$ for the braided vector subspace  of $\mathcal{K} $ given by 
the braided adjoint action of $\toba(V)$ on $\lambda_{a,b}$, that is
$$L(\lambda_{a,b}):=\operatorname{ad}_c \toba(V)(\lambda_{a,b})  = \big(\operatorname{ad}_c  \ku[x]/(x^{n})\big)(\Bbbk\{y\}) .$$

By definition, it holds that 
$L(\lambda_{a,b}) $ is linearly spanned by elements constructed from $y$ and the braided adjoint action of $x$. 
\begin{prop} \ 
\begin{enumerate}
 \item[$(a)$] $L(\lambda_{a,b})$ 
is a simple object of $\,_{T_{n,m}}^{T_{n,m}}\mathcal{YD}$, for all $0\leq a,b\leq nm-1$.
\item[$(b)$] The map $\lambda_{a,b}\mapsto L(\lambda_{a,b})$ yields
a bijective correspondence between the simple objects of  $\,_{\ku C_{nm}}^{\ku C_{nm}}\mathcal{YD}$ 
and the simple objects of $\,_{T_{n,m}}^{T_{n,m}}\mathcal{YD}$. 
\end{enumerate}
\end{prop}

\pf Follows from \cite[Proposition 3.5]{AHS} and \cite[Proposition 2.9]{AA}, respectively.
\epf

By the proposition above, we obtain a parametrization of the simple objects in $_{T_{n,m}}^{T_{n,m}}\mathcal{YD}$ and a 
way to compute them. The following result gives us the correspondence, through the braided equivalences described in Section \ref{sec:D-modules}, 
between these simple objects and the
parametrization obtained in 
Proposition \ref{prop:double-tnm}.

\begin{prop} \label{prop:corresp-L-V}
Let $0\leq i,j\leq nm-1$. Then 
$L(\lambda_{-i,-j}) \simeq F(\V)$.
 \end{prop}

\pf We prove the claim by showing explicitly the isomorphism. In order to construct the map, we describe first a 
distinguished basis of $L(\lambda_{a,b})$ for 
$\lambda_{a,b}=\Bbbk \{y\}$. 

Set $x_0=y$ and 
$x_{i+1}:=\operatorname{ad}_c x (x_i)=x\cdot x_i$ for all $i\geq 0$. By induction, one can prove that the 
following equality holds for all $k\geq 1$:
$$
x_{k}=x x_{k-1} -\omega^{k-1} \xi^{b} x_{k-1} x
=\sum_{\ell=0}^k(-1)^\ell  \xi^{\ell b} \omega^{\ell(\ell-1)/2} \binom{k}{\ell}_{\omega} x^{k-\ell}y x^{\ell}.
$$
In particular, this imples that $x_k=0$ for all $k\geq n$, because $w^n=1$ and $x^n=0$. 
Let $r$ the smallest positive integer such that $x_{r-1}\notin \Bbbk\{ x_0,\ldots,x_{r-2} \}$ and 
$x_{r}\in \Bbbk\{ x_0,\ldots,x_{r-1} \}$. Then, by the very definition of the elements $x_{k}$ 
and the module $L(\lambda_{a,b})$, the set
$\{x_0,\ldots,x_{r-1}\}$ is a linear basis of $L(\lambda_{a,b})$ and $x_{r}=0$. 
Now consider the following change of basis  given by 
\begin{align}\label{eq:muda-base}
z_0=x_0=y\quad\text{ and }\quad z_k=-\xi^{-b}\omega^{-(k-1)}\, x \cdot z_{k-1},\quad\text{ for all } 1\leq k\leq r-1.
\end{align}
On this basis, 
the action and coaction of $T_{n,m}$ on $L(\lambda_{a,b})$ is given by
\begin{align}\label{eq-action1}
	&g\cdot z_k=\xi^{b}\omega^k z_{k},\qquad\quad 
	 x\cdot z_k =
	\begin{cases}
\begin{array}{ll}
 - \xi^{b}\omega^k z_{k+1}& \text{ if } \quad 0\leq k<r-1;\\
 0& \text{ if } \quad k=r-1;
\end{array}
	\end{cases}\\
\label{eq-comulti}  
	&\rho (z_k)=\sum_{\ell=0}^{k}\beta'_{\ell,k}\,\,   x^\ell g^{a+k-\ell} \otimes z_{k-\ell},
	\end{align} 
for all $0\leq k \leq r-1$, where
$\beta'_{\ell,k}={k \choose \ell}_{\omega}\,\, \prod_{s= k-\ell}^{k-1} \left(\omega^{a}-\xi^{-b}\omega^{-s}\right)$ if $\ell\geq 1$ and $\beta'_{0,k}=1$.
%
%
Indeed, since $g\cdot x=\omega x$ and $x_k=x x_{k-1} -\omega^{k-1} \xi^{b} x_{k-1} x$, 
we obtain by induction on $k$ that  $g\cdot x_k=\xi^b\omega^kx_k$
for all $k\geq 0$. Consequently, $g\cdot z_k=\xi^b\omega^kz_k$ for all $0\leq k \leq r-1$ as claimed. 
Concerning the action of $x$, by \eqref{eq:muda-base} we have that 
$x\cdot z_{k}=(-\xi^{-b}\omega^{-k})^{-1}z_{k+1}=-\xi^{b}\omega^{k}z_{k+1}$ for all $0\leq k\leq r-2$.
If $k=r-1$, the claim follows since $x\cdot x_{r-1}=x_{r}=0$.
The assertion about the coaction follows by induction  
using that 
$\rho(z_{k+1})=-\xi^{-b}\omega^{-k}\rho(x\cdot z_k)$ and the compatibility condition.

\smallbreak
\noindent
{\bf Claim:} $r= \dim L(\lambda_{a,b}) =n$ if $m\nmid b$, and  $r=-a-\frac{b}{m} +1 \mod n$ if $m\mid b$. 

We prove first that if $r<n$, then necessarily $m\mid b$.
As $\dim L(\lambda_{a,b}) = r$, from \eqref{eq-comulti} we have that $0= \rho(x\cdot z_{r-1})$. But in such a case, the 
compatibility condition yields
\begin{align*}
 0 &= \rho(x\cdot z_{r-1}) = \sum_{\ell=0}^{r-1}\beta'_{\ell,r-1}\,\,   x_{(1)} x^\ell g^{a+r-1-\ell}\cS(x_{(3)}) \otimes x_{(2)}\cdot z_{r-1-\ell}\\
 &= \sum_{\ell=0}^{r-1}\beta'_{\ell,r-1}\,\,   x^{\ell+1} g^{a+r-1-\ell}\otimes z_{r-1-\ell} +
 \sum_{\ell=0}^{r-1}\beta'_{\ell,r-1}\,\,   g x^\ell g^{a+r-1-\ell}\otimes x\cdot z_{r-1-\ell}+\\
 & \quad + \sum_{\ell=0}^{r-1}\beta'_{\ell,r-1}\,\,   g x^\ell g^{a+r-1-\ell}\cS(x) \otimes g\cdot z_{r-1-\ell}\\
 &= \sum_{\ell=0}^{r-1}\beta'_{\ell,r-1}\,\,   x^{\ell+1} g^{a+r-1-\ell}\otimes z_{r-1-\ell} 
- \sum_{\ell=1}^{r-1}\beta'_{\ell,r-1}\,\,   g x^\ell g^{a+r-1-\ell}\otimes \xi^{b}\omega^{r-1-\ell}  z_{r-\ell}+\\
 & \quad + \sum_{\ell=0}^{r-1}\beta'_{\ell,r-1}\,\,   g x^\ell g^{a+r-1-\ell}(-g^{-1}x) \otimes \xi^{b}\omega^{r-1-\ell}z_{r-1-\ell}\\ 
 &= \beta'_{r-1,r-1}\,\,   x^{r} g^{a}\otimes z_{0} + \sum_{\ell=1}^{r-1}\beta'_{\ell-1,r-1}  x^{\ell} g^{a+r-\ell}\otimes z_{r-\ell} 
- \sum_{\ell=1}^{r-1}\beta'_{\ell,r-1}\xi^{b}\omega^{r-1}    x^\ell g^{a+r-\ell}\otimes z_{r-\ell}+\\
 & \quad - \sum_{\ell=0}^{r-2}\beta'_{\ell,r-1}\,\,   \omega^{a+r-1}x^{\ell+1} g^{a+r-1-\ell} \otimes \xi^{b}\omega^{r-1-\ell}z_{r-1-\ell}
 -\beta'_{r-1,r-1}\,\,   \omega^{a+r-1} \xi^{b} x^{r} g^{a} \otimes z_{0}
\end{align*}
\begin{align*}
  &=\beta'_{r-1,r-1}(1-\omega^{a+r-1} \xi^{b}) x^{r} g^{a} \otimes z_{0} +\\
  & \qquad +  \sum_{\ell=1}^{r-1}\big(\beta'_{\ell-1,r-1}-\beta'_{\ell,r-1}\xi^{b}\omega^{r-1} - 
  \beta'_{\ell-1,r-1}\xi^{b}\omega^{r-\ell}  \omega^{a+r-1}\big)  x^{\ell} g^{a+r-\ell}\otimes z_{r-\ell}.
\end{align*}
Thus $\beta'_{r-1,r-1}(1-\omega^{a+r-1} \xi^{b})=0$ and 
$\beta'_{\ell-1,r-1}-\beta'_{\ell,r-1}\xi^{b}\omega^{r-1} - 
  \beta'_{\ell-1,r-1}\xi^{b}\omega^{r-\ell} \omega^{a+r-1} =0$ for all $1\leq \ell \leq r-1$. 
 If $\beta'_{r-1,r-1} \neq 0$, then $1=\omega^{a+r-1} \xi^{b}$ and consequently $b+m(a+r-1) \equiv 0 \mod nm$ and 
 $m\mid b$ as desired. 
 On the other hand, $0=\beta'_{r-1,r-1} =\prod_{s= 0}^{r-2} \left(\omega^{a+s}\xi^{b}-1\right)\xi^{-b}\omega^{-s}$ 
 if and only if there exists $0\leq s\leq r-1$ such that 
  $1=\xi^{b}\omega^{a+s}$. As before, this implies that
 $m\mid b$. In particular, $\dim L(\lambda_{a,b}) = n$ if $m\nmid b$.
 
 Now we prove that $r=-a-\frac{b}{m} +1 \mod n$ if $m\mid b$. 
 Since for 
  all $1\leq \ell \leq r-1$ it holds that
  $\beta'_{\ell,r-1} = \frac{(r-\ell)_{\omega}}{(\ell)_{\omega}}(\omega^{a}-\xi^{-b}\omega^{-(r-1-\ell)})\beta'_{\ell-1,r-1}$, 
 from the equation above we get that   
  $$
  0 = \beta'_{\ell-1,r-1}\big(1-\omega^{\frac{b}{m}+r-1}  \frac{(r-\ell)_{\omega}}{(\ell)_{\omega}}(\omega^{a}-\omega^{-\frac{b}{m}-(r-1-\ell)})
  - \omega^{\frac{b}{m}+ r-\ell} \omega^{a+r-1}\big).  
  $$
Since $1= \beta'_{0,r-1}$, we necessarily have that 
\begin{align*}
& &0 &=1- \omega^{\frac{b}{m}+r-1}  \frac{(r-1)_{\omega}}{(1)	_{\omega}}\omega^{a}+
\frac{(r-1)_{\omega}}{(1)_{\omega}}\omega- \omega^{\frac{b}{m}+r-1}\omega^{a+r-1} \\
\Leftrightarrow & &  
1+ (r-1)_{\omega}\, \omega& =  
\omega^{\frac{b}{m}+r-1}  (r-1)_{\omega}\, \omega^{a}+
 \omega^{\frac{b}{m}+r-1}\omega^{a+r-1} \\
\Leftrightarrow & &  
1+ (r-1)_{\omega}\, \omega& =  \omega^{\frac{b}{m}+ a+r-1} \left( (r-1)_{\omega}+
 \omega^{r-1}\right) \\
 \Leftrightarrow & &  
(r)_{\omega}& =  
\omega^{\frac{b}{m}+a+r-1}  (r)_{\omega}.
 \end{align*}
 If $r<n$, then $(r)_{\omega}\neq 0$, and the equality above holds 
 if and only if $1=\omega^{\frac{b}{m}+a+r-1}$.  
 This implies that $r=-a-\frac{b}{m} +1 \mod n$ as desired. If $r=n$, a similar argument shows that $n=-a-\frac{b}{m} +1 \mod n$.

Finally, consider the linear map $\varphi: F(\V) \to L(\lambda_{-i,-j})$ given by $\varphi(v_{k})= z_{k}$ for all $0\leq k\leq r-1$.
From \eqref{double-tnm1}, \eqref{double-tnm2}, \eqref{eq-action1} and \eqref{eq-comulti}, it follows that 
$\varphi$ is an isomorphism in $_{T_{n,m}}^{T_{n,m}}\mathcal{YD}$.
\epf

\subsection{Finite-dimensional Nichols algebras over the simple modules $\V$} 
In this subsection we isolate the pairs $(i,j)$ with $0\leq i,j\leq nm-1$ such that $\toba(\V)$ is finite-dimensional.
Recall that, by Proposition \ref{prop:corresp-L-V}, a simple module $\V \in \ydhnm$ corresponds to the simple module
$L(\lambda_{-i,-j}) \in\ _{T_{n,m}}^{T_{n,m}}\mathcal{YD}$ via a braided equivalence, and that $T_{n,m} \simeq \toba(V)\#\Bbbk C_{nm}$.
The next result is a direct consequence of 
\cite[Theorem 1.1]{AA} and Proposition \ref{prop:corresp-L-V}.

\begin{teor}\label{teo:nichols-finite-dimension} 
For all pairs $(i,j)$ with $0\leq i,j\leq nm-1$, it holds that
$\dim \toba(\V)<\infty$ if and only if $\dim \toba(V\oplus \lambda_{-i,-j})<\infty$.
\qed
\end{teor}

Here, $W_{-i,-j} =V\oplus \lambda_{-i,-j}  \in  \,_{C_{nm}}^{C_{nm}}\mathcal{YD}$ is a braided vector space of diagonal type 
associated with the basis $\{x,y\}$. 
Thus, one may decide whether $\toba(V\oplus \lambda_{-i,-j})$ is finite-dimensional or not by looking at the classification 
provided by Heckenberger \cite{He, H2}, which establishes the correspondence between finite-dimensional Nichols algebras of
diagonal type and certain generalized Dynking diagrams. In our case, this reads as follows.

The braiding $c$ associated with 
$W_{-i,-j}$ is given by:
\begin{align*}
c(x\otimes x)&=\xi^{m}x\otimes x,& c(x\otimes y)&=\xi^{-j}y\otimes x, \\ 
c(y\otimes x)&=\xi^{-im}x\otimes y,& c(y\otimes y)&=\xi^{ij}y\otimes y.
\end{align*}
Hence, the corresponding braiding matrix and generalized Dynkin diagram $D_{i,j}$ are  
\begin{align}\label{diag:dynkin}
\mathbf{q} = 
\left(\begin{smallmatrix}
 \xi^{m} & \xi^{-j}\\
 \xi^{-im} & \xi^{ij}
\end{smallmatrix}\right) & & 
\begin{picture}(40,10)(20,25)
\put(-20,23){$D_{i,j}:$}
\put(10,20){$\circ$}
\put(10,30){\tiny{$\xi^{m}$}}
\put(14.7,22.7){\line(1,0){50}} 
\put(29.7,30){\tiny{$\xi^{-j-im}$}}
\put(64.7,20){$\circ$} 
\put(64.7,30){\tiny{$\xi^{ij}$}}
\end{picture}
& 
\end{align}

By Theorem \ref{teo:nichols-finite-dimension} above, and since the braided vector spaces
$W_{-i,-j}$ are of rank 2 for all $i,j$, it is enough to settle when the generalized 
Dynkin diagram in \eqref{diag:dynkin} belongs to \cite[Table 1]{H2}.  We do this for some specific cases 
and in an schematic way. It might be worth noting that this procedure can be done for all families of semisimple 
modules in $\ydhnm$, eventually looking at braided vector spaces of higher rank.

The generalized Dynkin diagrams corresponding to arithmetic root systems of rank 1 and 2  
are listed in \cite[Table 1]{H2} by rows and depend on fixed parameters 
$q,\zeta \in \ku^{\times}$. We call them \textit{Heckenberger diagrams} for short.
In each row, we order the diagrams from left to right, and  
we set $\HK_{k,\ell}$ for the $\ell$-th diagram in the $k$-th row. Following Heckenberger's notation,
we index the first row by 0; it corresponds to a generalized Dynking diagram of rank 1. 
For instance, $\HK_{4,1}$ denotes the diagram  
\begin{picture}(70,20)(5,20)
\put(10,20){$\circ$}
\put(10,30){\tiny{$q$}}
\put(14.7,22.7){\line(1,0){50}} 
\put(32.7,30){\tiny{$q^{-2}$}}
\put(64.7,20){$\circ$} 
\put(64.7,30){\tiny{$-1$}}
\end{picture},
%
%
with $q\in \Bbbk^{\times}\smallsetminus \{1,-1\}$ and $q\in \bG'_4$. The set of indices of Heckenberger 
diagrams is the following:
\begin{align*}
 \mathcal{I} & := \{(0,1),(1,1), (2,1),(2,2), (3,1), (4,1), (4,2), (5,1), (5,2), (6,1),(6,2), (7,1),(7,2),  \\
& \qquad  (7,3), (7,4),(7,5),(8,1), (8,2), (8,3), (9,1), (9,2),(9,3),(10,1), (11,1), (11,2),  \\
& \qquad (11,3), (12,1), (12,2), (12,3), (12,4),(13,1),(13,2),(14,1),(14,2),(14,3), (14,4),\\
& \qquad (15,1), (15,2), (15,3),(15,4), (16,1),(16,2)\}.
 \end{align*}

\bigbreak
\begin{center}
 \textit{From now on we assume that $n=2$ and $m\geq 2$.}
\end{center}
\bigbreak

As $n=2$, the generalized Dynkin diagram $D_{i,j}$ in \eqref{diag:dynkin} reduces to the following diagrams
\begin{equation}\label{diag:dynkin-n=2}
\begin{aligned}
\begin{picture}(40,20)(20,25)
\put(-80,23){$D_{i,j}:$}
\put(-30,20){$\circ$}
\put(-30,30){\tiny{$-1$}}
\put(24.7,20){$\circ$} 
\put(24.7,30){\tiny{$\xi^{ij}$}}
\end{picture}
&  \text{ or }&
\begin{picture}(40,20)(20,25)
\put(40,23){$D_{i,j}:$}
\put(90,20){$\circ$}
\put(90,30){\tiny{$-1$}}
\put(94.7,22.7){\line(1,0){50}} 
\put(110.7,30){\tiny{$\xi^{-j-im}$}}
\put(144.7,20){$\circ$} 
\put(144.7,30){\tiny{$\xi^{ij}$}}
\end{picture}
\end{aligned}
 \end{equation}
\bigbreak

Besides, by Definition \ref{def:simplemodulesD}, we have that $\dim \V \leq 2$ for all $0\leq i,j \leq 2m-1$. In particular,
$\dim \V= 1$ if and only if $j=0$ and $i$ is even, or $j=m$ and $i$ is odd.

We start by setting some notations and conventions. 
Recall that $\xi\in\Bbbk^{\times}$ is a fixed primitive $2m$-root of unity; in particular $\xi^{m}= -1$.


\begin{lema} \label{lema:edge}
If $\xi^{-j-im}=-1$, then $\D_{i,j}$ is not a Heckenberger diagram.
\end{lema}
\pf If $\xi^{-j-im}=-1$, then $m\equiv j+im \, \mo\, 2m$. Thus, $mi\equiv ji+mi^2\equiv ji +mi \, \mo\, 2m$, which implies
that $ij\equiv 0 \, \mo\, 2m$. In particular, the label of one vertex of $D_{i,j}$ equals 1. The claim follows since
no Heckenberger diagram has a vertex labelled by $1$.
\epf

\begin{lema} \label{lema:rank1}
If $\xi^{-j-im}=1$, then $m$ divides $j$ and 
$\D_{i,j}$ is the union of  two Heckenberger diagrams of rank one. Moreover, $ \dim \V = i+\frac{j}{m} + 1 \mod n$, and
$\dim \toba(\V) <\infty$ if and only if $ n\nmid i\frac{j}{m}$. 
\end{lema}

\pf If $\xi^{-j-im}=1$ then $j\equiv -im \mod 2m$, which implies that $m$ divides $j$; say $j=mt$ for some
$t\in \Z$ and $\xi^{ij} = (\xi^{m})^{ti}$. Then, by Definition \ref{def:simplemodulesD} we have that 
$ \dim \V = i+t + 1 \mod n$.
Under this assumptions, it holds that 
$c_{V,\lambda_{-i,-j}}c_{\lambda_{-i,-j},V} = \id$, which implies by \cite[Theorem 2.2]{Gr} that 
$\toba(V\oplus \lambda_{-i,-j}) \simeq \toba(V)\otimes \toba(\lambda_{-i,-j})$ 
 as $\N_{0}$-graded braided objects. 
 Write $d=  n/\gcd (n,it)$ for the order of $(\xi^{m})^{it}$. Then, $\toba(\lambda_{-i,-j})$
 is finite-dimensional if and only if $d\neq 1$. In such a case,
 $\toba(\lambda_{-i,-j}) \simeq \Bbbk[y]/(y^{d}) $. 
  Since $\toba(V) \simeq \Bbbk[x]/(x^{n})$, 
 the lemma follows by Theorem \ref{teo:nichols-finite-dimension}.
\epf

\bigbreak
\begin{obs}\label{rmk:rank1}
Assume that $\xi^{-j-im}=1$. Then, by Theorem \ref{teo:nichols-finite-dimension} and 
Lemma \ref{lema:rank1} we have that
 $\toba(\V)$ is finite-dimensional 
if and only if $j=m$  and $i$ is odd. In such a case, $\xi^{ij}= \xi^{-mi^{2}}= (-1)^{i^{2}} =-1 $ and
the corresponding Heckenberger diagram is disjoint and equals
\begin{picture}(50,20)(5,20)
\put(10,20){$\circ$}
\put(10,30){\tiny{$-1$}}
\put(44.7,20){$\circ$} 
\put(44.7,30){\tiny{$-1$}}
\end{picture}. 
Set $\mathcal{V}_{i,m} = \Bbbk v_{0}$.
By the proof of Proposition \ref{prop:double-tnm}, we have that the 
module and comodule structure
of $\mathcal{V}_{i,m}$ in $\ _{H_{2,m}}^{H_{2,m}}\mathcal{YD}$ is given by $A\cdot v_{0} = \xi^{i} v_{0}$, $X\cdot v_{0} = 0$
and $\rho(v_{0}) = A^{m}\ot v_{0}$. Thus, the corresponding braiding 
 is $c(v_{0}\ot v_{0}) = A^{m}\cdot v_{0}\ot v_{0} = \xi^{im}  v_{0}\ot v_{0} = - v_{0}\ot v_{0}$.
 This implies that $\toba(\mathcal{V}_{i,m}) \simeq \Bbbk[x]/(x^{2})$ for all $i$ odd.
 \end{obs}
 
From now on, we  assume now that $\xi^{-j-im}\neq 1$, that is, 
the corresponding Heckenberger diagram in \eqref{diag:dynkin-n=2} is connected.
Consider the following subset of indices of $\mathcal{I}$:
 \begin{align*}
 \mathcal{L} & := \{(2,1),(2,2),(4,1), (4,2), (6,1),(6,2),(7,2), (7,3), (7,4),(7,5),(9,2),(9,3), (11,2), \\
& \qquad (11,3), (12,3), (12,4),(13,1),(13,2),(14,1),(14,2),(14,3),(14,4), (15,3),(15,4)\}.
 \end{align*}

 \begin{lema}\label{lem:hecken-pairs} 
 Let $0\leq i,j\leq 2m-1$ and assume $D_{i,j}$ is a connected Heckenberger diagram $\HK_{k,\ell}$. Then 
$(k,\ell) \in \mathcal{L}$.
\end{lema}
\pf  We show that the remaining cases are not possible.
Since a vertex of $D_{i,j}$ is labelled by $-1$, the cases in the subset
$\{(7,1),(8,1),(9,1),(11,1),(12,1),(12,2),(15,1),(15,2)\}$ are not possible.
Besides, the cases in $\{(1,1),(3,1),(5,1),(5,2),(10,1)\}$ are neither possible by Lemma \ref{lema:edge}.

Suppose that $\D_{i,j}= \HK_{16,1}:$
\begin{picture}(70,10)(5,20)
\put(10,20){$\circ$}
\put(10,30){\tiny{$-1$}}
\put(14.7,22.7){\line(1,0){50}} 
\put(32.7,30){\tiny{$-\zeta^{-3}$}}
\put(64.7,20){$\circ$} 
\put(64.7,30){\tiny{$-\zeta$}}
\end{picture},
with $\zeta\in \bG'_{7}$. Then $\xi^{ij} = -\zeta$ and $\xi^{-j-im} = -\zeta^{-3}$. This implies that 
$14 ij\equiv 0\, \mo\, 2m$, $7ij\not \equiv 0\, \mo\, 2m$, $14 j\equiv 0 \, \mo\, 2m$ and $-7mi -7j\not \equiv 0  \, \mo\, 2m$. In particular,
$7j=m s$ for some $s\in \Z$. Thus $0\not\equiv 7ij\equiv msi\, \mo\, 2m$, which implies that $i$ and $s$ are odd. 
But on the other hand, $-7mi -7j=-7mi-ms=m(-7i+s)\equiv 0\,\mo\, 2m$, a contradiction. 
The cases $\HK_{16,2}$,  $\HK_{8,2}$ and $\HK_{8,3}$ can be treated similarly and are left to the reader.
\epf 

\begin{obs}
\label{rmk:dual} 
By direct inspection on \cite[Table 1]{H2}, one sees that 
there exists $(k,\ell_1)$ such that $\D_{i,j}=\HK_{k,\ell_1}$ if and only if there exists $\ell_2$ such that 
$\D_{-(i+1),-(j+m)}=\HK_{k,\ell_2}$. In such a case, it holds that $\ell_1\neq \ell_2$.
\end{obs}

In the following proposition we state the converse of Lemma \ref{lem:hecken-pairs}. In particular, we 
present the list of all pairs $(i,j)$ such that $\toba(\V)$ is finite-dimensional, and the conditions 
that the integers $m, i$ and $j$ must satisfy.

\begin{prop}\label{prop:classifica}
For all pairs $(k,\ell)\in \mathcal{L}$ there exist $m$, $i$ and $j$ with 
$m\geq 2$ and $0\leq i,j\leq 2m-1$ such that $D_{i,j}= H_{k,\ell}$. In particular,
for any such triple $(m, i,j)$, it holds that $\dim\toba(\V)<\infty$. The list of all possible triples 
$(m,i,j)$ is displayed in Table \ref{tab:Heck} below. 
\hspace*{-0.5cm}
\def\arraystretch{1.4}
\begin{table}[ht]
	\caption{Heckenberger diagrams of rank 2 associated with $\V$. \newline
	Write $m_{2} = \gcd(m,j) $,
$m=m_{1}m_{2}$ with $m_{1}$, $m_{2}$ and $a,b \in \N$.}\label{tab:Heck}
\begin{tabular}{|c|c|c|c|}
	\hline
$m$  & $(i,j)\,:\,\D_{i,j}=\HK_{k,\ell}$ & Conditions & $(k,\ell)$ \\ \hline
$ m_1 m_2$  & $(1-m_1a,m_2(m_1-b) )$           &  $\gcd(m_1, b)=1$,$\,\,\,a,b\in \bI$, $\,\,\,m_1\neq 1$. & $(2,1)$\\ \hline
$ m_1 m_2$   &$(m_1a,m_2b)$                     &$\gcd(m_1,b)=1$, $\,\,a,b\in \bI$,$\,\,\,m_1\neq 1$.& $(2,2)$ \\ \hline
$ m_1 m_2$  & \multirow{2}{*}{$(\frac{m_1 a+1}{2} ,m_2 b)$} &$\gcd(m_1, b)=1$, $\,\,a,m_1\in \bI$, & $(4,1)$ \\
      &                            &$\,m_1\neq 1$,$\,\,\,b\,\equiv\frac{m_1 a+1}{2}\,\mo\, 2$. & $(4,2)$  \\ \hline
$3m_2$  & \multirow{2}{*}{$(3a-1, \, m_2 b)$} & \multirow{2}{*}{$\!\!\quad\gcd(3, b)=1$, $\,\,\,b\,\equiv a\,\mo\, 2$.}&$(6,1)$ \\
&                                     &   &$(6,2)$ \\ \hline
\multirow{2}{*}{$6 m_2$}& \multirow{2}{*}{$(4+12a,  \,  m_2b)$} & \multirow{2}{*}{$\quad \gcd(6, b)=1$.} &$(7,2)$ \\
 &                                     &  & $(7,3)$ \\ \hline
\multirow{2}{*}{$6 m_2$} & \multirow{2}{*}{$(9+12a,  \, m_2b)$}&\multirow{2}{*}{ $\quad \gcd(6, b)=1$.}&$(7,4)$  \\ 
&                                     & &$(7,5)$ \\ \hline
$9 m_2$ & $(12+18a,  \, m_2b)$                & $\,\,\,\gcd(3, b)=1$,$\,\,\,b\in \bP$. &$(9,2)$\\ \hline
$9 m_2$& $(7+18a, \, m_2 b)$               &$\,\,\,\gcd(3, b)=1$,$\,\,\,b\in \bI$. &$(9,3)$ \\ \hline
 $4 m_2$ & $ (2+8a, \, m_2 b)$               & $\quad \gcd(2, b)=1$. &$(11,2)$ \\ \hline
 $4 m_2$ & $(7+8a,  \,  m_2 b)$              & $\quad \gcd(2, b)=1$. &$(11,3)$\\ \hline
$12m_2$& $(8+24a, \, m_2 b)$             & $\!\!\quad \gcd(12,b)=1$. &$(12,3)$ \\ \hline
$12m_2$ & $(17+24a, \, m_2 b)$             & $\!\!\quad \gcd(12,b)=1$. &$(12,4)$\\ \hline
$5m_2$&$ (2+10a, \, m_2 b)$             & $\quad \gcd(5,b)=1$, $\quad b\in \bP$. &$(13,1)$ \\ \hline
$5m_2$ &$ (9+10a, \, m_2 b)$             & $\quad \gcd(5,b)=1$, $\quad b\in \bI$. &$(13,2)$ \\ \hline
\multirow{2}{*}{$10 m_2$} & \multirow{2}{*}{$(17+20a,\,m_2 b)$} & \multirow{2}{*}{ $\quad \gcd(10,b)=1$.}&$(14,1)$  \\
&                                     &  &$(14,2)$  \\ \hline
\multirow{2}{*}{$10 m_2$} & \multirow{2}{*}{$(4+20a,\,m_2 b)$} & \multirow{2}{*}{$\quad \gcd(10,b)=1$.}&$(14,3)$  \\
&                                     &   & $(14,4)$ \\\hline
\multirow{2}{*}{$15 m_2$} &$ (25+30a,\, m_2 b)$   \text{ or }   &$\quad \gcd(15,b)=\gcd(15,c)=1$, & \multirow{2}{*}{$(15,3)$}\\
                           &$ (10+30a,\,m_2c)$                  &$b\in \bP$ \text{ and } $c\in \bI$. &\\ \hline
\multirow{2}{*}{$15 m_2$}&$ (6+30a, \,  m_2 b)$ \text{ or }     &$\quad \gcd(15,b)=\gcd(15,c)=1$, &  \multirow{2}{*}{$(15,4)$} \\ 
                    &$(21+30a,\,m_2c)$                       &$b\in \bI$ \text{ and } $c\in \bP$. &\\ \hline
\end{tabular}
\end{table}
\end{prop}

\pf We prove only the first $4$ cases, the remaining ones follow \textit{mutatis mutandis}. 

\bigbreak
\noindent $\diamond\,\, $ {\bf Case $(2,2)$:} 
$\HK_{2,2}=$
\begin{picture}(70,10)(5,20)
\put(10,20){$\circ$}
\put(10,30){\tiny{$-1$}}
\put(14.7,22.7){\line(1,0){50}} 
\put(36.7,30){\tiny{$q$}}
\put(64.7,20){$\circ$} 
\put(64.7,30){\tiny{$-1$}}
\end{picture},
with $q\in \ku^{\times}$, $q\neq \pm 1$. 
Let $(i,j)$ be such that $\D_{i,j}=\HK_{2,2}$. 
Set $m_2:=\gcd(m,j)$ and write $j=m_2b$, $m=m_1m_2$ with $\gcd(b,m_1)=1$ and $m_1,b \in \Z$.
Then $(m,i,j) = (m_{1}m_{2}, m_{1}a, m_{2}b)$, with $m_{1}\neq 1$ and $a,b\in \bI$. Indeed, 
since $\xi^{ij}=-1$, we must have that  $ij\equiv m\,\mo\, 2m$, say 
$ij=m+2mk =m(2k +1)$ for some $k\in \Z$; in particular $ib=m_1(2k+1)$. As $\gcd(m_1,b)=1$, 
it follows that $i=m_1a$ for some $a\in \Z$ and whence $m_1ab=m_1(2k+1)$. 
Then $ab=2k+1$ and consequently $a,b\in \bI$. Suppose that $m_1=1$. Then $m=m_{2}$, $j=mb$ and 
$\xi^{j+im}=\xi^{m(b+i)}=\pm 1$ which is a contradiction.
Conversely, any triple $(m, i,j)= (m_{1}m_{2}, m_{1}a, m_{2}b)$ with the assumptions above satisfies that $\D_{i,j}=\HK_{2,2}$:
clearly $\xi^{ij}=-1$ and $\xi^{j+im}= q \neq \pm 1$, since $\xi^{j+im}= \pm 1$ if and only if $m$ divides $m_2b$, 
and the latter occurs if and only if 
$m_1=1$.

\bigbreak
\noindent $\diamond\,\, $ {\bf Case $(2,1)$:} 
$\HK_{2,1}=$
\begin{picture}(70,10)(5,20)
\put(10,20){$\circ$}
\put(10,30){\tiny{$q$}}
\put(14.7,22.7){\line(1,0){50}} 
\put(36.7,30){\tiny{$q^{-1}$}}
\put(64.7,20){$\circ$} 
\put(64.7,30){\tiny{$-1$}}
\end{picture},
where $q\in \ku^{\times}$, $q\neq \pm 1$.
By Remark \ref{rmk:dual}, it follows that $\D_{i,j}=\HK_{2,1}$ 
if and only if $\D_{1-i,m-j}=\HK_{2,2}$. By the previous case, we must have that 
$m = m_{1}m_{2}$,
$1-i \equiv m_{1}a \mod 2m$ and $m-j \equiv m_{2}b \mod 2m$ with $m_{1}\neq 1$ and $a,b\in \bI$. Thus,
$(m,i,j) = (m_{1}m_{2}, 1-m_{1}a, m-m_{2}b)$ with $m_{1}\neq 1$ and $a,b\in \bI$.

\bigbreak
\noindent $\diamond\,\, $ {\bf Case $(4,1)$:} 
$\HK_{4,1}=$
\begin{picture}(70,10)(5,20)
\put(10,20){$\circ$}
\put(10,30){\tiny{$q$}}
\put(14.7,22.7){\line(1,0){50}} 
\put(36.7,30){\tiny{$q^{-2}$}}
\put(64.7,20){$\circ$} 
\put(64.7,30){\tiny{$-1$}}
\end{picture},
where $q\in \ku^{\times}$, $q\neq \pm 1$ and $q\notin \bG_4$.
Let $(i,j)$ be such that $\D_{i,j}=\HK_{4,1}$. As above, 
set $m_2:=\gcd(m,j)$ and write $j=m_2b$, $m=m_1m_2$ with $\gcd(b,m_1)=1$ and $m_1,b \in \Z$.
Then
$(m,i,j) = ( m_{1}m_{2}, \frac{m_{1}a+1}{2}, m_{2}b)$ with 
$m_{1}$, $a \in \bI$ and $b\equiv  \frac{m_{1}a+1}{2}\,\mo\, 2$.    
Since $q=\xi^{ij}$ and $q^{-2} = \xi^{j+im}$, we have that $ij\not\equiv 0,m \mo 2m$ and
$\xi^{j+im}=\xi^{-2ij}$. So, $-2ij\equiv j+im \,\mo\, 2m$ and whence $-(2i+1)j=(2k-i)m$, for some $k\in \Z$. 
Thus, $-(2i+1)b=(2k-i)m_1$ and $m_1$ divides $-(2i+1)$. Let $a\in \Z$ such that 
$-(2i+1)=m_1a$. Then $j\neq 0 \neq b$, $m_1,a$ are odd integers, $i=-(m_1a+1)/2$ and $ab=(2k-i)$. 
Hence $b\equiv i\,\mo\, 2$. The converse is straightforward.

\bigbreak
\noindent $\diamond\,\, $ {\bf Case $(4,2)$:} 
$\HK_{4,2}=$
\begin{picture}(70,10)(5,20)
\put(10,20){$\circ$}
\put(10,30){\tiny{$-q^{-1}$}}
\put(14.7,22.7){\line(1,0){50}} 
\put(36.7,30){\tiny{$q^{2}$}}
\put(64.7,20){$\circ$} 
\put(64.7,30){\tiny{$-1$}}
\end{picture},
where $q\in \ku^{\times}$, $q\neq \pm 1$ and $q\notin \bG_4$.
This case follows from the previous one by replacing $q$ by $-q^{-1}$. 
\epf 

\subsection{Presentations}
\label{subsec:presentationNichols}

In this subsection we present the finite-dimensional Nichols algebras corresponding to the first 
five rows of Table \ref{tab:Heck} by generators and relations. The braidings are described in the proofs.
We freely use the results and notation
from \S \ref{subsec:Nichols} and \cite{AGi}, where braidings of rank 2 were described.
We also provide the PBW-bases, the dimensions and the Hilbert series; for the latter, we write 
$(n)_{t} = t^{n-1} +\cdots + t + 1$ for $t \in \Bbbk$ and $n\in \N$. 

In general, the description of finite-dimensional Nichols algebras by generators and relations is a (quite) hard problem. 
For small dimensions, one can use the help of \texttt{GAP}, but for the general case one would need to use more theoretical tools to solve the 
problem, such as convex orders or 
Weyl groupoids as in the abelian case, see \cite{An}. We computed the first cases by using theoretic arguments and lots of computations.
For space (and time) reasons, we leave the description of the remaining cases to future work.

To give the presentation, as in the proof of Proposition \ref{prop:classifica},
we fix the following conventions: for $0\leq i,j\leq 2m-1$ set 

\bigbreak
\begin{center}
$m_2:=\gcd(m,j)$, $m=m_1m_2$ with $\gcd(b,m_1)=1$ and $m_1,b \in \Z$. 
\end{center}
\bigbreak

The following description corresponds to the first row of Table \ref{tab:Heck}. 
\begin{prop}\label{prop:Nichols2.1} 
Assume $(m,i,j) = (m_{1}m_{2}, 1-m_{1}a, m_{2}(m_{1}-b))$ with $m_{1}\neq 1$ and 
$a,b \in \bI$.
\begin{enumerate}
 \item If $m_1$ is odd then 
	\begin{align*}
	\toba(\V)\simeq \ku\langle v_0,v_1\,:\,v_0^{m_1}=0,\,v_0v_1-v_1v_0=0,\,v_1^2-(1+\xi^{-m_2b})^2v_0^2=0\rangle.
	\end{align*}
\item If $m_1$ is even then 
\begin{align*}
\toba(\V)\simeq \ku\langle v_0,v_1\,:\,v_0^{2m_1}=0,\,v_0v_1+v_1v_0=0,\,v_1^2+(1-\xi^{-m_2b})^2v_0^2=0\rangle.
\end{align*}
\end{enumerate}
In particular, for both cases we have that 
$ \{v^{a_0}_0v^{a_1}_1 : \,  0\leq a_0 < \gcd(2,m_1)\, m_1, \, 0\leq a_1\leq 1\} $ is a PBW-basis of 
$ \toba(\V) $, $ \dim\toba(\V) = \gcd(2,m_1)\, 2m_1 $ and the Hilbert series is $ (\gcd(2,m_1)\, m_1 )_t \, (2)_t $.
	
\end{prop}
\pf  If $m_1$ is odd, then the braiding of $\V$ is given by
\begin{align*}
(c(v_i \ot v_j))_{i, j \in \{0,1\} } = \begin{pmatrix}
-\xi^{-m_2b} v_0\otimes v_0 & v_1\otimes v_0 - (\xi^{-m_2b}+1) v_0\otimes v_1 \\
-\xi^{-m_2b}v_0\otimes v_1 & - v_1\otimes v_1+(1-\xi^{-2m_2b})(\xi^{-m_2b}+1) v_0\otimes v_0
\end{pmatrix};
\end{align*}
this is a braiding of type $\hit_{1, 2}$.
Similarly, when $m_1$ is even the braiding reads
\begin{align*}
(c(v_i \ot v_j))_{i, j \in \{0,1\} } = \begin{pmatrix}
\xi^{-m_2b} v_0\otimes v_0 & -v_1\otimes v_0 + (\xi^{-m_2b}-1) v_0\otimes v_1 \\
-\xi^{-m_2b}v_0\otimes v_1 & - v_1\otimes v_1+(1-\xi^{-2m_2b})(\xi^{-m_2b}-1) v_0\otimes v_0
\end{pmatrix},
\end{align*}
which is a braiding of type $\hit_{1, 2} \, (a) \, (b)$.
Hence, the corresponding Nichols algebras, PBW-bases and dimensions 
are given by \cite[Proposition 3.10 and 3.11]{AGi}, respectively. 
\epf

The next presentation corresponds to the second row of Table \ref{tab:Heck}. Its proof 
also follows from \cite[Proposition 3.10 and 3.11]{AGi}.

\begin{prop}\label{prop:Nichols2.2}
Assume $(m,i,j) = (m_{1}m_{2}, m_{1}a, m_{2}b)$ with $m_{1}\neq 1$ and 
$a,b \in \bI$. 
\begin{enumerate}
 \item  If $m_1$ is odd then  
		\begin{align*}
		\toba(\V)\simeq \ku\langle v_0,v_1\,:\,v_0^2=0,\,\,v_1^{m_1}=0,\,v_0v_1+\xi^{-m_2b}v_1v_0=0\rangle.
		\end{align*}	
\item If $m_1$ is even then 
\begin{align*}
\toba(\V)\simeq \ku\langle v_0,v_1\,:\,v_0^2=0,\,\,v_1^{2m_1}=0,\,v_0v_1+\xi^{-m_2b}v_1v_0=0\rangle.
\end{align*}
\end{enumerate}
In particular, in both cases $ \{v^{a_0}_0v^{a_1}_1 :  \, 0\leq a_0\leq 1, \,  0\leq a_1 < \gcd(2,m_1)\, m_1\} $ is a PBW-basis of 
$ \toba(\V) $, $ \dim\toba(\V) = \gcd(2,m_1)\, 2m_1 $ and the Hilbert series is $ ( \gcd(2,m_1)\, m_1 )_t \, (2)_t $.
\qed
\end{prop}

We split the description of the Nichols algebras corresponding to 
the third row into two cases. 
Assume $(m,i,j) = (m_{1}m_{2}, \frac{m_{1}a+1}{2}, m_{2}b)$ with $a,m_{1}\in \bI$, $m_{1}\neq 1$ and 
$b \equiv\frac{m_1 a+1}{2} \mo\, 2$.

\begin{prop}\label{prop:Nichols4.1_beven} 
Suppose that $b$ is even.
Set 
$ p = \xi^{\frac{mab+m_2b}{2}}$, $k = (1-\xi^{2m_2b})\xi^{\frac{mab-m_2b}{2}}(1-\xi^{m_2b})$ and write $ v_{10} = [v_1, v_0]_c =  v_1v_0 - m\,  c(v_1\ot v_0)  = v_1v_0 - p\, v_0v_1$.
Then $ \toba(\V)$ is the algebra generated by $v_0, v_1$ satisfying the following relations: 
\begin{align} \label{eqn:relation1_row3_beven}
&v_{10}v_0 = v_0v_{10},\qquad \qquad v_1v_{10} = -p^2\, v_{10}v_1 -kp(p-1)v_0^3, \\\label{eqn:relation3_row3_beven}
&v_0^{m_1} = 0, \\\label{eqn:relation4_row3_beven}
&v_1^{2m_1} = 0 .
\end{align}
In particular, $ \{v^{a}_0v^{b}_{10}v^{c}_1 :  \, 0\leq a < m_1, \,  0\leq b < 2, \,  0\leq c < 2m_1\} $ is a PBW-basis of 
$ \toba(\V) $, $ \dim\toba(\V) = 4m_1^2 $ and the Hilbert series is $ ( m_1 )_t \, (2)_{t^2} \, (2m_1 )_t  $.

\end{prop}
\pf 
The braiding associated with $\V$ is of type $\hit_{1, 2}$ and equals
\begin{align*}
(c(v_i \ot v_j))_{i, j \in \{0,1\} } = \begin{pmatrix}
p\,  v_0\otimes v_0 & p^{-1} \, v_1\otimes v_0 + (p-p^{-1})v_0\otimes v_1 \\
p\, v_0\otimes v_1 & -p^{-1}\,  v_1\otimes v_1+k\, v_0\otimes v_0
\end{pmatrix}.
\end{align*}
 
As $ p \neq \pm 1 $, by \cite[Proposition 3.10]{AGi} it follows that $ \toba(\V) $ does not have any quadratic relation.
By computing $\Ker Q_{3}$,  we obtain that 
relations $ \eqref{eqn:relation1_row3_beven} $  must hold 
in $\toba(\V)$. Also, 
$ \eqref{eqn:relation3_row3_beven} $ holds since $v_{0}^{m_{1}}$ 
is a primitive element of degree bigger than $1$. To show that 
\eqref{eqn:relation4_row3_beven} holds we use skew-derivations. Write $\partial_{0} = \partial_{v_{0}^{*}}$
and $\partial_{1} = \partial_{v_{1}^{*}}$ for the skew-derivations associated with the dual basis $\{v_{0}^{*},v_{1}^{*}\}$
of $\V^{*}$.

Consider the pre-Nichols algebra given by the quotient $ \widetilde{\toba} (\V) = T(\V)/I $, where $ I $ 
is the two-sided ideal generated by the relations \eqref{eqn:relation1_row3_beven} and  \eqref{eqn:relation3_row3_beven}.
In particular, in $ \widetilde{\toba} (\V)  $ it holds that 
$
	 v_{10}^2 =  \dfrac{kp(p-1)^2}{p^2+1}\, v_{0}^4
$.

We claim that  $ B= \{v_0^{a}\,  v_{10}^{b}\, v_1^{c}\,:\,  0\leq a\leq m_1-1,  \ 0\leq b\leq 1, 0\leq c\} $  
spanns linearly $ \widetilde{\toba} (\V)  $. 
Indeed, since $v_1\, v_0^{a} = (a)_p \, v_0^{a-1} v_{10} +p^a \, v_0^{a} v_1$ for all $0\leq a \leq m_1-1$, we have that
\begin{align*} 
 v_1\, v_0^{a} \,  v_{10}^{b}\, v_1^{c} & = (a)_p \, v_0^{a-1} v_{10} v_{10}^b v_{1}^c + 
 p^a \, v_0^{a} v_1  v_{10}^b v_{1}^c \\
& =\begin{cases}
\begin{array}{cc}
(a)_p \, v_0^{a-1} v_{10} v_{1}^c + p^a \, v_0^{a} v_{1}^{c+1} & \text{if }b=0, \\
\frac{kp(1-p)(p^{a+2}+1)}{p^2+1} \, v_0^{a+3}  v_{1}^c - p^{a+2} \, v_0^{a} v_{10} v_{1}^{c+1} & \text{if }b=1.
\end{array} 
\end{cases}
\end{align*}
This implies that the linear subspace spanned by $ B $ is a left ideal of $ \widetilde{\toba} (\V) $. 
Since $ \widetilde{\toba} (\V) $ is 
generated by $v_{0}$ and $v_{1}$, the claim follows.

\smallbreak
Set $ q = -p^{-1} $ and for $ \ell\in \mathbb{N}, \iota\in \mathbb{Z}_2 $
consider $ \Lambda_\iota^{(\ell)}\in \widetilde{\toba} (\V)$ defined recursively by 
$$\begin{cases}
\begin{array}{lcl}
\Lambda_\iota^{(\ell)}=a^{(\ell)}_{\iota,0} \,\, v_0\Lambda_{\iota}^{(\ell-1)}+a^{(\ell)}_{\iota,1}\,\, 
v_1 \Lambda_{1-\iota}^{(\ell-1)},  &&  
	\iota\in \mathbb{Z}_2,\, \ell \geq 2, \\
	\Lambda_\iota^{(1)} = v_\iota, &&  \iota\in \mathbb{Z}_2,
\end{array}
 \end{cases}$$
 where
\begin{align*}
	a^{(\ell)}_{\iota,0} =\begin{cases}
	\begin{array}{cc}
	k^{1-\iota} \, q^{\ell+\iota-2}  & \text{if $\ell$ is even}, \\
	k^{\iota} \, q^{\ell-\iota-1}  & \text{if $\ell$ is odd},
	\end{array} 
	\end{cases} \text{ and } \qquad a^{(\ell)}_{\iota,1} = (q-p)^{1-\iota}(-q)^{\iota+\ell-2}.
\end{align*}
By induction, one may prove that the evaluation of $\Lambda_\iota^{(\ell)}$ on $\partial_{0}$ and $\partial_{1}$ equals
\begin{align*}
	\partial_0(\Lambda_\iota^{(\ell)})= a^{(\ell)}_{\iota,0} \,\, (\ell)_q \, \Lambda_{\iota}^{(\ell-1)} \qquad\text{ and }\qquad 
	\partial_1(\Lambda_\iota^{(\ell)})= a^{(\ell)}_{\iota,1} \,\, (\ell)_q  \Lambda_{1-\iota}^{(\ell-1)} 
	\qquad \forall\ \iota\in\mathbb{Z}_2, \,\ell\geq 2.
\end{align*}
As $ q\in \bG'_{2m_1} $, it follows that $ \Lambda_0^{(2m_1)}\in\J(\V)$. On the other hand, 
if we write $ \Lambda_0^{(2\ell)} $ in terms of the 
set $ B $, we get that
\begin{align*}
	\Lambda_0^{(2\ell)} &= v_1^{2\ell} (q^2+1)^\ell q^{\ell(2\ell-3)} \\
	 &\quad  - \sum_{r=1}^{\ell} v_0^{2r}v_1^{2(\ell-r)} (q^2+1)^{\ell-r} (-k)^r q^{2(\ell-r)(\ell -1-r) -\ell+1}\binom{2\ell}{2r}_q (2r-1)_q \\
	&  \quad + \sum_{r=1}^{\ell -1} v_0^{2r-1}v_{10}v_1^{2(\ell-r)-1} (q^2+1)^{\ell-r} (-k)^r q^{2(\ell -1-r)^2 -\ell+1} \binom{2\ell}{2r+1}_q \dfrac{(2r)_q}{q+1}.
\end{align*}
In particular, $ \Lambda_0^{(2m_1)} = v_1^{2m_1} (q^2+1)^{m_1} q^{m_1(2m_1-3)} -  v_0^{2m_1}  (-k)^{m_1} q^{-m_1+1} (2m_1-1)_q = 
v_1^{2m_1} (q^2+1)^{m_1} q^{m_1(2m_1-3)}$, which implies that $v_1^{2m_1} \in \mathcal{J}(\V)$.

Finally, let $ J $ be the two-sided ideal of $T(\V)$ generated by the relations 
$ \eqref{eqn:relation1_row3_beven} -\eqref{eqn:relation4_row3_beven} $. 
Clearly, $T(\V)/J \simeq  \widetilde{\toba} (\V) /(v_{0}^{2m_{1}})$ is
a pre-Nichols algebra of $\V$ linearly spanned by the  
set $ B'= \{v_0^{a}\,  v_{10}^{b}\, v_1^{c}\,:\,  0\leq a\leq m_1-1,  \ 0\leq b\leq 1, 0\leq c\leq 2m_{1}-1\} $.
To show that $ \toba (\V) \simeq T(\V)/J $, it suffices to see that for all $x\in T(\V)/J $ it holds that 
$\partial_{0}(x) = 0 = \partial_{1}(x)$ if and only if $x=0$. This follows from the equalities
\begin{align*}
\partial_1(v_0^{a}\,  v_{10}^{b}\, v_1^{c}) &= (c)_q (-q^2)^b(-q)^a v_0^{a}\,  v_{10}^{b}\, v_1^{c-1}, \\ 
\partial_0(v_0^{a}\,  v_{10}^{b}) &= (a)_p v_0^{a-1}\,  v_{10}^{b} + p^a(1-p^2)\delta_{b,1} \, v_0^a v_1, 
\end{align*}
for all $ 0\leq a\leq m_1-1,  \ 0\leq b\leq 1, 0\leq c \leq 2m_1 -1 $. 
In particular, $ B' $ is a basis of $ \toba (\V)  $ and the related Hilbert series is $ ( m_1 )_t \, (2)_{t^2} \, (2m_1 )_t  $.
\epf

Now we describe the finite-dimensional Nichols algebra corresponding to the third row for the case when $b$ is odd.
The proof is quite similar to the one above.
\begin{prop}\label{prop:Nichols4.1_bodd} 
	Suppose that $b$ is odd. Set 
		$ p = \xi^{\frac{mab-m_2b}{2}}$, 
		$q=-p^{-1}$
		and $k = (1-\xi^{2m_2b})\xi^{\frac{mab-m_2b}{2}}(1+\xi^{m_2b})$.
		Write $ v_{10} = v_1v_0 + q\, v_0v_1$.
	Then $ \toba(\V)$ is the algebra generated by $v_0, v_1$ satisfying the following relations: 
	\begin{align} \label{eqn:relation1_row3_bodd}
		& v_{10}v_0 = -v_0v_{10}, \qquad \qquad v_1v_{10} = q^2\, v_{10}v_1 -kq(q-1)v_0^3,\\\label{eqn:relation3_row3_bodd}
		&v_0^{2m_1} = 0, \\\label{eqn:relation4_row3_bodd}
		&v_1^{m_1} + \dfrac{(-k)^{\frac{m_1-1}{2}} p^{\frac{1-m_1^2}{2}}}{(p^2+1)^{\frac{m_1-1}{2}}(p+1)} v_0^{m_1-2}v_{10} = 0.
	\end{align}	
	In particular, $ \{v^{a}_0v^{b}_{10}v^{c}_1 :  \, 0\leq a < 2m_1, \,  0\leq b < 2, \,  0\leq c < m_1\} $ is a PBW-basis of 
	$ \toba(\V) $, $ \dim\toba(\V) = 4m_1^2 $ and the related Hilbert series is $ ( 2m_1 )_t \, (2)_{t^2} \, ( m_1 )_t  $.
\end{prop}

\pf First note that 
the braiding associated with $\V$ is of type $\hit_{1, 2}(a) \, (b)$ (see \cite{AGi}):
\begin{align*}
(c(v_i \ot v_j))_{i, j \in \{0,1\} } = \begin{pmatrix}
q\,  v_0\otimes v_0 & p \, v_1\otimes v_0 + (p+q)v_0\otimes v_1 \\
-q\, v_0\otimes v_1 & p\,  v_1\otimes v_1+k\, v_0\otimes v_0
\end{pmatrix}.
\end{align*}

As $ p \neq \pm 1 $, by \cite[Proposition 3.11]{AGi} we have that $ \toba(\V) $ does not have any quadratic relation.
Relations $ \eqref{eqn:relation1_row3_bodd} $ follow from a direct
computation of $\Ker Q_{3}$, and 
$ \eqref{eqn:relation3_row3_bodd} $ follows from the fact that $v_{0}^{2m_{1}}$ is primitive, as the order of $q$ is $2m_{1}$.

Consider now the pre-Nichols algebra given by the quotient $ \widetilde{\toba} (\V)= T(\V)/I $, where $ I $ is the 
two-sided ideal generated by the relations 
\eqref{eqn:relation1_row3_bodd} and  \eqref{eqn:relation3_row3_bodd}.
In particular, in $ \widetilde{\toba} (\V)  $ it holds that 
$v_{10}^2 =  -\dfrac{kq(q-1)^2}{q^2+1}\, v_{0}^4$.

Moreover, following the proof of  Proposition \ref{prop:Nichols4.1_beven}, one can prove that 
$ B= \{v_0^{a}\,  v_{10}^{b}\, v_1^{c}\,:\,  0\leq a\leq 2m_1-1, \, 0\leq b\leq 1, \, 0\leq c\} $  spanns linearly $ \widetilde{\toba} (\V)  $. 

To show that \eqref{eqn:relation4_row3_bodd} holds in the Nichols algebra, we use skew-derivations on some fixed elements:
for $ \ell\in \mathbb{N}$ and $\iota\in \mathbb{Z}_2 $ consider the elements $ \Lambda_\iota^{(\ell)}\in \widetilde{\toba} (\V)$ defined recursively by
$$\begin{cases}
\begin{array}{lcl}
\Lambda_\iota^{(\ell)}=a^{(\ell)}_{\iota,0} \,\, v_0\Lambda_{\iota}^{(\ell-1)}+a^{(\ell)}_{\iota,1}\,\, v_1 \Lambda_{1-\iota}^{(\ell-1)},  
& & \iota\in \mathbb{Z}_2,\, \ell \geq 2, \\
\Lambda_\iota^{(1)} = v_\iota, & & \iota\in \mathbb{Z}_2,
\end{array}
\end{cases}
$$
where
\begin{align*}
a^{(\ell)}_{\iota,0} =\begin{cases}
\begin{array}{cc}
k^{1-\iota} \, p^{\ell+\iota-2}  & \text{if $\ell$ is even}, \\
k^{\iota} \, p^{\ell-\iota-1}  & \text{if $\ell$ is odd},
\end{array} 
\end{cases} \text{ and } \qquad a^{(\ell)}_{\iota,1} = (p-q)^{1-\iota}p^{\iota+\ell-2}.
\end{align*}
By induction, we may prove that for all $\iota\in\mathbb{Z}_2$, and  $\ell\geq 2$ one gets
\begin{align*}
\partial_0(\Lambda_\iota^{(\ell)})= a^{(\ell)}_{\iota,0} \,\, (\ell)_p \, \Lambda_{\iota}^{(\ell-1)}\qquad \text{ and }\qquad 
\partial_1(\Lambda_\iota^{(\ell)})= a^{(\ell)}_{\iota,1} \,\, (\ell)_p  \Lambda_{1-\iota}^{(\ell-1)},
\end{align*}
which implies that $ \Lambda_1^{(m_1)}\in\J(\V)$, since $ p\in \bG'_{m_1} $. 
Writing $ \Lambda_1^{(2\ell+1)} $ in terms of $ B $ yields
\begin{align*}
\Lambda_1^{(2\ell+1)} &= v_1^{2\ell +1} (p^2+1)^\ell p^{\ell(2\ell -1)} \\
& - \sum_{r=1}^{\ell} v_0^{2r}v_1^{2(\ell-r)+1} (p^2+1)^{\ell-r} (-k)^{r} p^{2(\ell -r)^2 -\ell+1} \binom{2\ell+1}{2r}_p (2r-1)_p \\
&   - \sum_{r=1}^{\ell } v_0^{2r-1}v_{10}v_1^{2(\ell-r)} (p^2+1)^{\ell-r} (-k)^{r} p^{2(\ell -r)(\ell -1-r) -\ell+2} \binom{2\ell+1}{2r+1}_p \dfrac{(2r)_p}{p+1}.
\end{align*}
In particular, $ \Lambda_1^{(m_1)} = v_1^{m_1} (p^2+1)^{\frac{m_1-1}{2}} p^{\frac{(m_1-1)(m_1-2)}{2}} + 
v_0^{m_1-2}v_{10} \dfrac{(-k)^{\frac{m_1-1}{2}} p^{\frac{m_1+3}{2}}}{p+1} $, from which relation $ \eqref{eqn:relation4_row3_bodd} $ follows.

Let $J$ be the two-sided ideal of $T(\V)$ generated by the relations 
\eqref{eqn:relation1_row3_bodd} --\eqref{eqn:relation4_row3_bodd}. Clearly, the pre--Nichols algebra 
$ T(\V)/J$ is isomorphic to $\widetilde{\toba} (\V) / J'$, where
$J'$ is the two-sided ideal generated by the relation \eqref{eqn:relation4_row3_bodd}, and it is linearly spanned by the image
of the set $B$.
To show that $T(\V)/J \simeq \toba(\V)$ it is enough 
to prove that for all $x\in T(\V)/J$ it holds that 
$\partial_{0}(x) = 0 = \partial_{1}(x)$ if and only if $x=0$. But this follows from the following equalities for 
$ 0\leq a\leq 2m_1-1$ and  $0\leq b\leq 1, 0\leq c \leq m_1 -1 $:
\begin{align*}
\partial_1(v_0^{a}\,  v_{10}^{b}\, v_1^{c}) &= (c)_p p^{a+2b} v_0^{a}\,  v_{10}^{b}\, v_1^{c-1}, \\ 
\partial_0(v_0^{a}\,  v_{10}^{b}) &= (a)_q v_0^{a-1}\,  v_{10}^{b} + q^a(q^2-1)\delta_{b,1} \, v_0^a v_1.
\end{align*} 
\epf

We present next the Nichols algebras corresponding to the fourth row.

\begin{prop}\label{prop:Nichols6.2} 
Assume that $(m,i,j) = (3m_{2},  3a-1, m_{2}b)$ with $\gcd(3,b)=1$ and $b\equiv a \mod 2$.
Set $p = (-1)^a\xi^{-m_2b}\in\bG_3'$, $q=p^{-1}$, 
$ k = (-1)^{a}(p-1)$ and $v_{10}=v_1v_0-(-1)^{a+1}p\,v_0v_1$.
Then $ \toba(\V)$ is the algebra generated by $v_0$ and $v_1$ satisfying:
\begin{align}
\label{eqn:relation1_row4}
  0&= v_0^3 ,\\
\label{eqn:relation2_row4}
   0&=v_{10}v_0+(-1)^{a}(1+2q)\, v_0v_{10}-2\,  v_0^2v_1,\\ 
 \label{eqn:relation3_row4}
   0&= v_1v_{10}+(-1)^{a}(2+q)\, v_{10}v_1-2 p\,  v_0v_1^2,\\
\label{eqn:relation4_row4}
   0&= v_1^3+(-1)^{a}(2q-2)\, v_0^2v_1+(1+2q) \, v_0v_{10}.   
\end{align}
In particular, $ B= \{v_0^{a_0}\,  v_{10}^{a_1}\, v_1^{a_2}\,:\,  0\leq a_0, a_2\leq 2,\ 0\leq a_1\leq 1\} $ 
is a PBW-basis of $\toba(\V)$, $ \dim\toba(\V) = 18 $ and the Hilbert series is $ (3)_t^2 \, (2)_{t^2} $.
\end{prop}

\pf The braiding associated  with $\V$ is given by
\begin{align*}
	(c(v_i \ot v_j))_{i, j \in \{0,1\} } = \begin{pmatrix}
		p\, v_0\otimes v_0 & \, (-1)^{a}q v_1\otimes v_0 + (p+q)\, v_0\otimes v_1 \\
		(-1)^{a+1}p\, v_0\otimes v_1 & q\, v_1\otimes v_1 + k\, v_0\otimes v_0
	\end{pmatrix}.
\end{align*}
If $ a $ is odd, then it is a braiding of type $\hit_{1, 2}$, and if $ a $ is even, then it is a braiding of type $\hit_{1, 2} \, (a) \, (b)$, see  \cite{AGi}. 
Clearly, the associated  Nichols algebra has no quadratic relation since $ p, q\in\bG_3' $.
A quick calculation shows that the cubic relations  \eqref{eqn:relation1_row4} -- \eqref{eqn:relation4_row4}
correspond to elements in $\Ker Q_{3}$.

Let $ \widetilde{\toba} (\V) = T(\V)/I $ be the pre-Nichols algebra 
given by the quotient of $T(\V)$ by the two-sided ideal
generated by the relations \eqref{eqn:relation1_row4} -- \eqref{eqn:relation4_row4}. 
We show that actually $ \widetilde{\toba} (\V) \simeq  \toba (\V)$. 
First note that, from these relations it follows that the equality $v_{10}^2 = 2q\, v_0v_{10} v_1$  holds in $ \widetilde{\toba} (\V)$.

Let $ B= \{v_0^{a_0}\,  v_{10}^{a_1}\, v_1^{a_2}\,:\,  0\leq a_0, a_2\leq 2,\ 0\leq a_1\leq 1\} $. Then, 
by showing that any element of the form $ v_1 \, v_0^{a_0}\, v_{10}^{a_1}\, v_1^{a_2} $ belongs to the linear spann of $B $, 
a straightforward computation similar to the one in the proof of Proposition \ref{prop:Nichols4.1_beven} yields that 
$B$ spanns linearly $ \widetilde{\toba} (\V) $. 

As a consequence, we get that 
the highest degree of a (non-zero) element in $ \widetilde{\toba} (\V) $ is $ 6 $.
Then, by the Poincar\'e duality, all possible relations appear in degree 2 or 3.
As we computed all quadratic and cubic relations, we have that 
$ \toba (\V) \simeq  \widetilde{\toba} (\V) $ and $ B $ is a basis of $ \toba (\V) $. Consequently,
the related Hilbert series is $ ( 3 )_t^2 \, (2)_{t^2} $.
\epf

We end this subsection with the presentation of the finite-dimensional 
Nichols algebras corresponding to the fifth row.
\begin{prop}\label{prop:Nichols7.2} 
Assume that $(m,i,j) = (6m_{2},  4+12a, m_{2}b)$ with $\gcd(6,b)=1$ and $ a \in \Z$.
Set $p = \xi^{4m_2b}\in\bG_3'$, $q=-\xi^{3m_2b}\in\bG_4'$, $ k = p\, (p+q)$, $v_{10}=v_1v_{0}-p\,v_{0}v_1$ and $v_{110}=v_1v_{10}-pq\,v_{10}v_1$.
Then $ \toba(\V)$ is the algebra generated by $v_0$ and $v_1$ satisfying:
\begin{align}\label{eqn:relation1_row5}
		v_0^3 &= 0, \\\label{eqn:relation2_row5}
		v_{10}v_0 &= -p(q+1)\, v_0v_{10} - p^2q\, v_0^2v_1, \\ \label{eqn:relation3_row5} 
		 v_1v_{110} &= (1-p-pq)\, v_{110}v_{1}  +(p-1)(1-pq)\, v_{10}v_1^2\\
		 \notag &\ \ \ \, + p(p-1)\, v_0v_1^3  +p(pq-1)\, v_0^2v_{10}, \\ \label{eqn:relation4_row5}
		 v_1^4 &= -(p+q)^2((2p^2-p+q(2-p))\, v_{0}^2v_{1}^2  +(p-1-q)\, v_0v_{10}v_1 \\
		 \notag &\hspace*{2,3cm} - (1+pq)\, v_{10}^2 +(1-p-p^2q)\, v_0v_{110}).
	\end{align}
	In particular, $ B= \{v_0^{a_0}\,  v_{10}^{a_1}\, v_{110}^{a_2}\, v_1^{a_3}\,:\,  0\leq a_0, a_1\leq 2,\ 0\leq a_2\leq 1,\ 0\leq a_3\leq 3\} $ 
	is a PBW-basis of $\toba(\V)$, $ \dim\toba(\V) = 72 $ and the Hilbert series is $ (4)_t\, (3)_t\, (3)_{t^2} \, (2)_{t^3} $. 
\end{prop}
\pf 
The braiding of $\V$ is of  type $\hit_{1, 2}$  and equals
\begin{align*}
(c(v_i \ot v_j))_{i, j \in \{0,1\} } = \begin{pmatrix}
p\, v_0\otimes v_0 & -q\, v_1\otimes v_0 + (p+q)\, v_0\otimes v_1 \\
p\, v_0\otimes v_1 & q\, v_1\otimes v_1 + k\, v_0\otimes v_0
\end{pmatrix}.
\end{align*}
A sheer calculation of $\Ker Q_{2}$, $\Ker Q_{3}$ and $\Ker Q_{4}$ give us the relations \eqref{eqn:relation1_row5} -- \eqref{eqn:relation4_row5}.
Thus, we need only to show that the pre-Nichols algebra given by the quotient 
$ \widetilde{\toba} (\V) = T(\V)/I $, where $ I $ is the two-sided ideal generated by these relations, is actually the Nichols algebra. 
In particular, in $ \widetilde{\toba} (\V)  $ we have the following equalities
\begin{align*}
v_{110}v_{0}&=  (1+p)(1+q)\, v_0v_{110} -p(1+2q)\, v_{10}^2 -v_0v_{10}v_1 -p(2+q)\, v_0^2v_1^2, \\
v_{110}v_{10} &= (1-pq+q) v_{10}v_{110} +(-2p-2q-1) v_{10}^2v_1 + (p+q+pq)v_0v_{110}v_1 \\ 
& \ \ \ +(1+2p^2q)v_0v_{10}v_1^2 + (2p-q+1)v_0^2v_1^3, \\
v_{10}^3 &= -pq \, v_0 v_{10} v_{110} - (p + q) v_0 v_{10}^2v_1 + (p-1) v_0^2 v_{110} v_1 + (1-p-pq) v_0^2 v_{10} v_1^2, \\ 
v_{110}^2 &= (p+2-3pq) v_{10}v_{110}v_1 + (3pq-2p-1) v_{10}^2v_1^2 + (p+q+2pq-1) v_0v_{110}v_1^2 \\
& \ \ \ + (1-p-5q-4pq) v_0v_{10}v_1^3 + (4pq-2q+6p^2)  v_0^2v_{10}^2.
\end{align*}
From these equations, one can see that the set  
$ B= \{v_0^{a_0}\,  v_{10}^{a_1}\, v_{110}^{a_2}\, v_1^{a_3}\,:\,  0\leq a_0, a_1\leq 2,\ 0\leq a_2\leq 1,\ 0\leq a_3\leq 3\} $
spanns linearly $ \widetilde{\toba} (\V) $. 
Thus the top degree of $ \widetilde{\toba} (\V) $ is $12$.
Using the Poincar\'e duality, we obtain that $ \toba (\V) \simeq  \widetilde{\toba} (\V) $ unless there exists
a new relation of degree less or equal than 6.
With the use of \texttt{GAP}, one may compute explicitly $\Ker Q_{5}$ and $\Ker Q_{6}$
and check that no new relation appears. In particular, $ B $ is a basis of $ \toba (\V)$.
\epf

\subsection{Summary for $n=2$ and $m=2,3$}
In this last subsection, we collect the presentations for the cases $n=2$ and $m=2,3$. In particular,
we obtain the description of the finite-dimensional Nichols algebras corresponding to indecomposable modules
over the Hopf algebras $H_{n,m}= R_{n,m}^{*}$, with $n$ and $m$ as above. This follows from our previous results and 
as a consequence of the following theorem due to Andruskiewitsch and Angiono.

\begin{teor}\label{thm:indec-simple}\cite[Theorem 1.2]{AA}
Let $Z \in\ \ydhnm$ be such that $\dim \toba(Z) < \infty$. Then $Z$ is semisimple.
\end{teor}

As a first application, we obtain another version of \cite[Theorem A]{GJG}, see also \cite{X}.
\begin{teor}\label{thm:22}
Let $H_{2,2}$ be the dual of the Radford algebra of dimension $8$ and $V \in\ _{H_{2,2}}^{H_{2,2}}\mathcal{YD}$ an indecomposable module.
Then $\dim\toba(V) <\infty$ if and only if $V$ is simple and isomorphic either to $\mathcal{V}_{1,2}$, $\mathcal{V}_{3,2}$, 
$\mathcal{V}_{3,1}$, $\mathcal{V}_{3,3}$, $\mathcal{V}_{2,1}$ or $\mathcal{V}_{2,3}$.
In such a case, the corresponding Nichols algebras are presented as follows:
\begin{enumerate}
 \item[$(a)$] $\toba(\mathcal{V}_{1,2})\simeq \Bbbk[x]/(x^{2})$ with $A\cdot x=\xi x$, $X\cdot x =0$ and $\rho(x)=A^{2}\ot x$, 
 $\dim \toba(\mathcal{V}_{1,2})= 2$;
 \item[$(b)$] $\toba(\mathcal{V}_{3,2})\simeq \Bbbk[x]/(x^{2})$ with $A\cdot x=\xi^{3} x$, $X\cdot x =0$ and $\rho(x)=A^{2}\ot x$;
\item[$(c)$] $\toba(\mathcal{V}_{3,1})\simeq \ku\langle v_0,v_1\,:\,v_0^{4}=0,\,v_0v_1+v_1v_0=0,\,v_1^2+(1+\xi)^2v_0^2=0\rangle$, $\dim \toba(\mathcal{V}_{3,1}) =8$;
 \item[$(d)$] $\toba(\mathcal{V}_{3,3})\simeq \ku\langle v_0,v_1\,:\,v_0^{4}=0,\,v_0v_1+v_1v_0=0,\,v_1^2+(1-\xi)^2v_0^2=0\rangle$;
 \item[$(e)$]$\toba(\mathcal{V}_{2,1})\simeq \ku\langle v_0,v_1\,:\,v_0^2=0,\,\,v_1^{4}=0,\,v_0v_1-\xi v_1v_0=0\rangle$, $\dim \toba(\mathcal{V}_{2,1}) =8$;
 \item[$(f)$]  $\toba(\mathcal{V}_{2,3})\simeq \ku\langle v_0,v_1\,:\,v_0^2=0,\,\,v_1^{4}=0,\,v_0v_1+\xi v_1v_0=0\rangle$.
 \end{enumerate}
\end{teor}

\pf Let $V\in\ _{H_{2,2}}^{H_{2,2}}\mathcal{YD}$ be indecomposable and assume that $\dim V< \infty$. 
Then by Theorem \ref{thm:indec-simple}, $V$ is simple and isomorphic 
to some $\V$ with $0\leq i,j\leq 3$, since 
it is indecomposable. 

As $n=2=m$, we have that $\dim \V\leq 2$.
If $\dim V = 1$, the assertion follows from 
Remark \ref{rmk:rank1}, with $V$ isomorphic to $\mathcal{V}_{1,2}$ or $\mathcal{V}_{3,2}$. If $\dim V=2$,
then the claim follows by inspection on Table \ref{tab:Heck}, where one needs to look
only at the first three rows. In all cases, $m=m_{1}=2$ and  $m_{2}=1$.
The conditions on the first row implies that $V\simeq \mathcal{V}_{3,1}$ or
$V\simeq \mathcal{V}_{3,3}$, whereas conditions on the second row implies that $V\simeq \mathcal{V}_{2,1}$ or
$V\simeq \mathcal{V}_{2,3}$. Conditions on the third row are not compatible with our assumptions on $m$. 
Finally, the presentations follow from Remark \ref{rmk:rank1} and Proposition \ref{prop:Nichols2.1} and \ref{prop:Nichols2.2}.
\epf

We end the paper by describing the case for $n=2$ and $m=3$. In particular, 
we obtain the results of Hu and Xiong \cite{HX} on finite-dimensional Nichols algebras over indecomposable objects 
and we answer \cite[Question 1]{HX} explicitly.
The proofs follows \textit{mutatis mutandis} from the proof above. The corresponding presentations are referred to the results in the previous
section.

\begin{teor}\label{thm:23}
Let $H_{2,3}$ be the dual of the Radford algebra of dimension $12$ and $V \in\ _{H_{2,3}}^{H_{2,3}}\mathcal{YD}$ an indecomposable module.
Then $\dim\toba(V) <\infty$ if and only if $V$ is simple and isomorphic either to $\mathcal{V}_{1,3}$, $\mathcal{V}_{3,3}$, 
$\mathcal{V}_{5,3}$, which are 1-dimensional, or to $\mathcal{V}_{4,2}$, $\mathcal{V}_{4,4}$, $\mathcal{V}_{3,1}$, $\mathcal{V}_{3,5}$,
$\mathcal{V}_{2,2}$, $\mathcal{V}_{2,4}$, $\mathcal{V}_{5,1}$, $\mathcal{V}_{5,5}$, 
$\mathcal{V}_{5,2}$, $\mathcal{V}_{5,4}$, $\mathcal{V}_{2,1}$ or $\mathcal{V}_{2,5}$.
In such a case, the corresponding Nichols algebras are presented in 
\begin{enumerate}
 \item[$(a)$] Remark \ref{rmk:rank1} for $\mathcal{V}_{1,3}$, $\mathcal{V}_{3,3}$, 
$\mathcal{V}_{5,3}$;
 \item[$(b)$] Proposition \ref{prop:Nichols2.1} ($m_{1}$ odd) for $\mathcal{V}_{4,2}$ and $\mathcal{V}_{4,4}$;
\item[$(c)$] Proposition \ref{prop:Nichols2.2} ($m_{1}$ odd) for $\mathcal{V}_{3,1}$ and $\mathcal{V}_{3,5}$;
 \item[$(d)$] Proposition \ref{prop:Nichols4.1_beven} for $\mathcal{V}_{2,2}$ and $\mathcal{V}_{2,4}$;
 \item[$(e)$] Proposition \ref{prop:Nichols4.1_bodd} for $\mathcal{V}_{5,1}$ and $\mathcal{V}_{5,5}$;
 \item[$(f)$]  Proposition \ref{prop:Nichols6.2} for $\mathcal{V}_{5,2}$, $\mathcal{V}_{5,4}$, $\mathcal{V}_{2,1}$ and $\mathcal{V}_{2,5}$.
  \end{enumerate} \qed
\end{teor}


\end{document}